\pdfoutput=1
\documentclass[letterpaper,11pt]{article}
\usepackage[margin=1in]{geometry}
\usepackage[
            CJKbookmarks=true,
            bookmarksnumbered=true,
            bookmarksopen=true,
            colorlinks=true,
            citecolor=red,
            linkcolor=blue,
            anchorcolor=red,
            urlcolor=blue
            ]{hyperref}

\usepackage[title]{appendix}
\usepackage{algorithm}
\usepackage{algpseudocode}
\floatname{algorithm}{Algorithm}

\usepackage{subcaption}

\usepackage{graphicx}
\usepackage{booktabs} 
\usepackage{natbib}
\usepackage{amsmath,amsthm,amsfonts,amssymb}
\usepackage{hyperref}
\usepackage{color}
\usepackage{enumitem}
\usepackage{bm}
\usepackage{multirow}
\usepackage{bbm}
\usepackage{subfiles}
\usepackage{xspace}
\usepackage{caption}
\usepackage{mathtools} 
\usepackage{tikz}

\newtheorem{theorem}{Theorem}
\newtheorem{example}[theorem]{Example}
\newtheorem{definition}[theorem]{Definition}
\newtheorem{assumption}[theorem]{Assumption}
\newtheorem{corollary}[theorem]{Corollary}
\newtheorem{lemma}[theorem]{Lemma}
\newtheorem{proposition}[theorem]{Proposition}
\newtheorem{remark}[theorem]{Remark}
\theoremstyle{definition}

\providecommand{\keywords}[1]
{
  \small	
  \textbf{\textit{Keywords---}} #1
}

\usetikzlibrary{intersections,calc,decorations.pathreplacing,through,arrows}
\usetikzlibrary{ decorations.markings,positioning}

\tikzset{global scale/.style={
    scale=#1,
    every node/.append style={scale=#1}
  }
}

\newcommand{\prob}{\mathbb{P}}
\newcommand{\Expect}{\mathbb{E}}
\newcommand{\real}{\mathbb{R}}
\newcommand{\KL}{\textnormal{KL}}
\newcommand{\TV}{\textnormal{TV}}
\newcommand{\diff}{\textnormal{d}}
\newcommand{\op}{\textnormal{op}}
\newcommand{\pth}[1]{\left( #1 \right)}

\newcommand{\sth}[1]{\left\{ #1 \right\}}
\newcommand{\Norm}[1]{\|#1\|}

\newcommand{\calB}{{\mathcal{B}}}

\newcommand{\calF}{{\mathcal{F}}}
\newcommand{\calG}{{\mathcal{G}}}
\newcommand{\calH}{{\mathcal{H}}}

\newcommand{\calM}{{\mathcal{M}}}

\newcommand{\calP}{{\mathcal{P}}}

\newcommand{\calV}{{\mathcal{V}}}

\DeclareMathOperator{\tr}{tr}
\DeclareMathOperator{\diag}{diag}
\DeclareMathOperator{\Span}{span}

\DeclareMathOperator{\Bern}{Bern}
\DeclareMathOperator{\rank}{rank}
\DeclareMathOperator{\Ind}{\textbf{Ind}}

\newcommand{\kru}{{\mathsf{Kru}}}
\newcommand{\iid}{i.i.d.}
\newcommand{\iiddistr}{{\stackrel{\text{i.i.d.}}{\sim}}}

\title{
Identifiability and Estimation in High-Dimensional \\ Nonparametric Latent Structure Models
}
\author{Yichen Lyu and Pengkun Yang
\thanks{Accepted for presentation at the Conference on Learning Theory (COLT) 2025.}
\thanks{Y. Lyu and P. Yang are with Department of Statistics and Data Science, Tsinghua University.  
P.~Yang is supported in part by the National Key R\&D Program of China 2024YFA1015800.
}}

\begin{document}

\maketitle

\begin{abstract}%
  This paper studies the problems of identifiability and estimation in high-dimensional nonparametric latent structure models. We introduce an identifiability theorem that generalizes existing conditions, establishing a unified framework applicable to diverse statistical settings. 
  Our results rigorously demonstrate how increased dimensionality, coupled with diversity in variables, inherently facilitates identifiability.
  For the estimation problem, we establish near-optimal minimax rate bounds for the high-dimensional nonparametric density estimation under latent structures with smooth marginals.
  Contrary to the conventional curse of dimensionality, our sample complexity scales only polynomially with the dimension.
  Additionally, we develop a perturbation theory for component recovery and propose a recovery procedure based on simultaneous diagonalization. 
\end{abstract}


\section{Introduction}
High-dimensional statistical models play a pivotal role in modern statistics and are widely applied across diverse research domains. 
A central challenge in such settings is the notorious \emph{curse of dimensionality}: as dimensionality grows, the volume of the space expands exponentially, rendering data increasingly sparse. 
Consequently, reliable inference typically requires sample sizes that grow prohibitively with dimension, posing fundamental limitations in practice.

These challenges are starkly evident in high-dimensional nonparametric density estimation, where the absence of structural assumptions leads to slow convergence rates and severe data inefficiency. 
Yet in practice, such as generative models, underlying distributions often possess inherent structure that constrains the function space of interest.
Exploiting such a structure can circumvent the curse of dimensionality, enabling tractable estimation even in high-dimensional regimes.

A compelling example arises when high-dimensional data is generated by populations with latent subgroups exhibiting \emph{conditional independence}. 
Such models are prevalent in applications spanning medical diagnosis \cite{hall_nonparametric_2003}, image recognition \cite{juan_use_2002,juan_bernoulli_2004}, chemical and physical sciences \cite{kasahara_non-parametric_2014}. 
See \cite{chauveau_semi-parametric_2015} for a review. 
In bivariate problems, the structure reduces to a low-rank representation of the data matrix. 
Mathematically, the data distribution is modeled as
\begin{equation}
   \label{eq:MM}
   \mu=\sum_{k=1}^m\pi_k (\mu_{k1}\times \mu_{k2}\times\cdots \times \mu_{kd}),
\end{equation}
where $\pi_k> 0$ for $k\in[m]\triangleq\{1,\dots,m\}$, $\sum_{k=1}^m\pi_k=1$, $\mu_{k}\triangleq \mu_{k1}\times \mu_{k2}\times\cdots \times \mu_{kd}$ is a product measure on $\real^d$.
In this paper, we assume the number of components $m\ge 2$ is known and fixed.  
Methods for estimating $m$ are discussed in \cite{kasahara_non-parametric_2014}.  

This paper studies the central theoretical question concerning the identifiability of such mixture models and the estimation problem from a sample of $n$ independent and identically distributed (\iid) observations from $\mu$.
The model is said to be \emph{identifiable} if no other model within the family yields the same data distribution. 
For mixture models, only the mixing measure $\sum_{k=1} ^m\pi_k\delta_{\mu_k}$ can be uniquely identified \cite{chen_optimal_1995,heinrich_strong_2018,wu_optimal_2020}, where $\delta$ denotes the Dirac measure, and thus the components can be identified only up to a global permutation.

Suppose each component probability measure $\mu_k\in \calP_d$ for some family $\calP_d$, a necessary condition to ensure identifiability is that $\calP_d$ is a nonconvex set.
The families of distributions from many parametric models, such as Gaussians, are nonconvex by definition, whose identifiability has been extensively investigated. 
In the absence of explicit parametric assumptions, nonparametric models are often adopted in practice.
However, nonparametric families such as H\"{o}lder-smooth densities are convex, and the mixture models are less studied. 
In model~\eqref{eq:MM}, each component belongs to the \emph{nonconvex} family of product measures. 
Formally, we define the identifiability of~\eqref{eq:MM} as follows.
\begin{definition}[Identifiability]
\label{identi}
   Let $\mu=\sum_{k=1}^m\pi_k (\mu_{k1}\times \cdots \times \mu_{kd})$. 
   We say $\mu$ is identifiable if $\tilde{\mu}=\sum_{k=1}^m\tilde \pi_k(\tilde\mu_{k1}\times \cdots \times \tilde\mu_{kd})=\mu$ implies that  there exists a permutation $\sigma:[m]\mapsto [m]$ such that $\pi_k=\tilde\pi_{\sigma(k)},\mu_{kj}=\tilde{\mu}_{\sigma(k)j}$ for all $k\in[m]$ and $j\in[d]$. 
\end{definition}

\subsection{Gaps in the Identifiability Conditions of Existing Literature}
\label{Motivation}
We begin by reviewing previous results on the identifiability conditions for model~\eqref{eq:MM}. \cite{teicher_identifiability_1967} was among the first to investigate this topic for the parametric case, establishing an equivalence between the identifiability of high-dimensional mixtures of product measures and the identifiability of one-dimensional mixtures with an unknown number of components.
For the nonparametric settings, \cite{hall_nonparametric_2003} made a pioneering contribution by addressing the identifiability for $m=2$. 
A cornerstone result is provided by \cite{allman_identifiability_2009} as stated below.
\begin{theorem}[Linear Independence Condition]
\label{linearind}
    Suppose $d\geq 3$ and $\mu$ can be expressed as \eqref{eq:MM}. If, for each $j\in[d]$, $\mu_{1j},\dots,\mu_{mj}$ are linearly independent, then $\mu$ is identifiable.
\end{theorem}
Theorem \ref{linearind} builds on an algebraic result by \cite{kruskal_three-way_1977}, who established the uniqueness of the canonical polyadic (CP) decomposition for three-way tensors. We refer to \cite{Kolda_09} for a comprehensive review of tensor decomposition.
The linear independence condition has since become a foundational assumption in many studies developing algorithms for model \eqref{eq:MM}. Notable examples include \cite{benaglia_em-like_2009,levine_maximum_2011,AGH+2014,zheng_nonparametric_2020,lu_nonparametric_2022}.

While the linear independence condition is widely adopted as a standard assumption in existing algorithms, the condition does not hold in numerous scenarios, as shown in the examples below.

  \begin{example}[Conditional \iid\ Model]
  \label{ex:CI}
    In \eqref{eq:MM}, for each $k \in[m]$, $\mu_{k1} = \dots = \mu_{kd} $. Hence,
    \begin{equation*} 
    \mu = \sum_{k=1}^m \pi_k \mu_{k1}^{\times d}. 
    \end{equation*}
    The linear independence condition fails when $\mu_{11},\dots,\mu_{m1}$ are linearly dependent.
 \end{example}
 \begin{example}[Bernoulli Mixture Model]
    \label{ex:BMM}
    The distribution of each $\mu_{kj}$ in \eqref{eq:MM} is given by a Bernoulli distribution:
    \begin{equation*}
       \mu_{kj}=\Bern(\alpha_{kj}).
   \end{equation*}
   The linear independence condition fails when $m \geq 3$.
 \end{example}
Both examples are special cases of \eqref{eq:MM} and are important topics of independent interest. 
 The conditional \iid\ model is closely related to learning mixing measures from group observations and the sparse Hausdorff problems, as discussed in \cite{rabani_learning_2014,li_learning_2015,gordon_sparse_2020,wu_optimal_2020,fan_efficient_2023}. The Bernoulli mixture model has been extensively studied by theoretical computer scientists \cite{feldman_learning_2008,pmlr-v134-gordon21a,gordon_identification_nodate} and finds applications in areas such as text learning, image recognition, and image generation \cite{juan_use_2002,juan_bernoulli_2004}.

Although Theorem \ref{linearind} does not apply to these examples, the recent progress shows that the models can be identified under certain conditions on the dimensionality and the \emph{diversity} along each variable. 
For instance, \cite{tahmasebi_identifiability_2018} showed that under certain separability conditions, the Bernoulli mixture model with $d\geq 2m-1$ is identifiable. They further generalized this result to the finite support case.
For the conditional \iid\ model, \cite{vandermeulen_operator_2019} showed that $\mu$ is identifiable when $d\geq 2m-1$. 
Remarkably, despite the failure of the linear independence condition, the threshold $d=2m-1$ emerges as a valid criterion for identifiability.
In Section \ref{sec_identi}, we bridge the gap by providing general identifiability conditions for model~\eqref{eq:MM} when the linear independence does not necessarily hold. 
\subsection{Related Work on the Estimation Problem}
We also study the estimation problem for model \eqref{eq:MM} given a finite sample. It is well known that in the nonparametric setting, density estimation suffers from the curse of dimensionality \cite{tsybakov_introduction_2009}. However, for model \eqref{eq:MM}, the latent structure from conditional independence substantially reduces model complexity: whereas a generic density estimation problem typically exhibits exponential rate on the dimension $d$, we show in Section \ref{sec:complexity} that the complexity of model \eqref{eq:MM} depends only polynomially on $d$.

For the estimation of components, we establish a perturbation analysis under quantitative assumptions. 
Specifically, given an error bound between $\mu$  and its estimate 
$\hat{\mu}$, we aim to derive quantitative error bounds between the component distributions 
$\mu_{kj}$
  and their corresponding estimates 
$\hat{\mu}_{kj}$. 
Prior work has established perturbation results in several special cases. For example, \cite{hall_nonparametric_2003} derived an asymptotic result for the two-component case; 
\cite{bhaskara_uniqueness_2013} gives quantitative rates in concrete cases;
\cite{vandermeulen_operator_2019} proposed a spectral method for the conditional i.i.d.\ model with consistency guarantees; and \cite{gordon_identification_nodate} obtained near-optimal bounds for the Bernoulli mixture model. These results suggest that the error in estimating the components is of the same order as the error in estimating the full model, which motivates the general perturbation theory developed in Section \ref{sec_roc}.

Algorithmic development under general identifiability conditions is another interesting question. 
Existing algorithms are broadly categorized into two types. 
The first is based on the nonparametric Expectation-Maximization (NPEM) algorithm \cite{benaglia_em-like_2009,benaglia_bandwidth_2011,levine_maximum_2011,chauveau_semi-parametric_2015}. 
While this iterative method is straightforward to implement, it lacks global convergence guarantees and is sensitive to the initial model. 
The second approach treats the model as a high-order tensor and applies algorithms from tensor decomposition. 
Recent works \cite{gordon_hadamard_2022,gordon_identification_nodate} successfully applied this framework to Bernoulli mixture models. 
While tensor-based algorithms benefit from a robust theoretical foundation, they are typically limited to discrete cases.

To address this gap, several recent works have adapted tensor methods to continuous settings. 
For example, \cite{bonhomme_estimating_2016} truncated the orthogonal basis in the $L^2$ space and applied tensor decomposition techniques, with the convergence rate depending on the precision of the truncation. \cite{zheng_nonparametric_2020} introduced a method for selecting a finite functional basis under the linear independence condition, which can be estimated using kernel density estimators. 
\cite{lu_nonparametric_2022} combines these approaches, thereby reducing the error rates. 
The linear independence condition remains crucial in many existing algorithms.

\subsection{Our Contributions}
Motivated by the theoretical gaps presented in the previous subsections, we study the identifiability and estimation problem of model~\eqref{eq:MM}.
Our main contributions are as follows:
\begin{itemize}
    \item \emph{A general, unified identifiability theorem.} 
    We propose an identifiability theorem in Section \ref{sec_identi} that unifies and extends all the previous identifiability conditions for model~\eqref{eq:MM}. 
    Notably, our result explains why high-dimensional variables with \emph{diversity} aid the identifiability. 
    \item \emph{Quantitative rates of convergence.}  
    We establish a perturbation theory in Section \ref{sec_roc} for estimating the components under an \emph{incoherence} condition. 
    Moreover, we derive near-optimal minimax risk bounds for high-dimensional nonparametric density estimation, where the sample complexity scales only polynomially with the dimension.
    \item \emph{A recovery algorithm under incoherence conditions.} 
    We develop a recovery algorithm for model \eqref{eq:MM} in Section \ref{sec_algo} that operates from an estimator of the joint density close to the true density. 
    Our algorithm successfully recovers the component densities relying only on incoherence rather than linear independence. 
\end{itemize}

\noindent\textbf{Notations}
 Let $[n]\triangleq \{1,2,\dots ,n\}$. 
 Let $\Delta^{n-1}\triangleq \{(x_1,\dots ,x_n)\in\real^n:x_i\ge 0, \sum_{i=1}^n x_i=1\}$ denote the $n$-simplex.
 For $\alpha\in\real$, the Dirac measure on the point $\alpha$ is defined as $\delta_{\alpha}$.
 The operator $\otimes$ denotes the Kronecker product for vectors and matrices, and the tensor product in general Hilbert spaces. 
 For $f,g\in\calH$, the angle between them is denoted as $\theta(f,g)\triangleq \cos^{-1}\frac{\langle f,g\rangle}{\|f\|_2\|g\|_2}$.
 For $f,g\in L^2(\real)$, the inner product is defined as $\langle f,g\rangle = \int f(x)g(x)\diff x$.
 For a finite rank linear operator $T$, denote the $i$-th largest singular value of $T$ by $\sigma_{i}(T)$.
 For two matrices $A=(a_{ij}),B=(b_{ij}) \in\real^{m\times n}$, the Hadamard product is denoted as $A\circ B = (a_{ij}b_{ij})_{i,j=1}^{m,n}\in\real^{m\times n}.$ 
 For two positive sequences $\{a_n\}$ and $\{b_n\}$, we write $a_n\lesssim b_n$ if $a_n\leq C b_n$ for a constant $C$, and $a_n\asymp b_n$ if $a_n\lesssim b_n$ and $b_n\lesssim a_n$, and we write $a_n\lesssim_{q} b_n, a_n\asymp_{q} b_n$ to emphasize that the $C$ depends on a parameter $q$.

\section{Model Identifiability without Linear Independence}
\label{sec_identi}
In this section, we establish the identifiability condition for model \eqref{eq:MM}. 
Without additional assumptions on the model, the joint measure $\mu$ is generally not identifiable. 
For instance, when $d=2$ and $\mu_{kj}$'s are discrete, model \eqref{eq:MM} reduces to the low rank decomposition of a matrix, which is well-known to be nonunique. 
Furthermore, for $d\ge 3$, additional variables are not helpful without diversity conditions: if $\mu_{k1}=\mu_{1}$ for all $k\in[m]$, the joint measure then becomes
\begin{equation*}
    \label{counter_identi}
    \mu=\mu_{1}\times \left(\sum_{k=1}^m\pi_k( \mu_{k2}\times\cdots \times \mu_{kd})\right).
\end{equation*}
Suppose $X=(X_1,\dots,X_d)\sim \mu$. Then $X_1$ is independent of $(X_2,\dots, X_d)$ and the model needs to be identified by the remaining $d-1$ variables. 
The following definition quantifies the diversity of a variable $X_j$ via the set of the conditional distributions $\{\mu_{kj}:k\in[m]\}$.

\begin{definition}[$\ell$-Independence]
\label{def:l_independence}
Let $(X_1,\dots,X_d)\sim \mu$ for $\mu$ in~\eqref{eq:MM}.
We say the $j$-th variable is $\ell$-independent if every subset of $\{\mu_{kj}\}_{k=1}^m$ of cardinality $\ell$ is linearly independent. 
Let 
\[
\textnormal{\textbf{Ind}}_\mu(j)
\triangleq\max\{\ell: \text{$j$-th variable is $\ell$-independent} \}
\]
For a subset $S\subseteq [d]$, define $\textnormal{\textbf{Ind}}_\mu(S)\triangleq\sum_{j\in S}\textnormal{\textbf{Ind}}_\mu(j)$, and let \[
    \tau_{\mu}(S) \triangleq \min\{m,\textnormal{\textbf{Ind}}_\mu(S)-|S|+1\}
    \]
    denote the total excess independence in $S$.
\end{definition}
Definition \ref{def:l_independence} is a generalization of Kruskal rank to probability measures. 
As special cases, $\Ind_\mu(j) = 1$ corresponds to identical components, where $\mu_{1j} = \dots = \mu_{mj}$, while $\Ind_\mu(j) = m$ corresponds to full linear independence.
Definition~\ref{def:l_independence} captures an intermediate notion between these two extremes.
Similar concepts can be found in \cite[Definition 4.1]{vandermeulen_generalized_2022}.
In particular, $\Ind_\mu(j)=2$ is equivalent to $\mu_{1j},\dots,\mu_{mj}$ are pairwise distinct---a property we formally define below as the separability condition.
\begin{definition}[Separability Condition]
\label{separa}
    Let $(X_1,\dots,X_d)\sim \mu$ for $\mu$ in~\eqref{eq:MM}. 
    The $j$-th variable is said to be separable if $\mu_{kj}\neq \mu_{k'j}$ for every pair of distinct indices $k\ne k'\in [m]$. 
    We denote by $N(\mu)$ the number of separable variables in model~\eqref{eq:MM}. 
\end{definition}
We now state our main result for the identifiability condition based on $\ell$-independence.
\begin{theorem}
   \label{identi_thm}
     Let $\mu$ be defined as in \eqref{eq:MM}. If there exists a partition $S_1,S_2,S_3$ of $[d]$ satisfying
    \begin{equation}
    \label{condition}
        \tau_{\mu}(S_1)+\tau_{\mu}(S_2)+\tau_{\mu}(S_3)\geq 2m+2,
    \end{equation}
        then $\mu$ is identifiable. 
        Conversely, there exists a non-identifiable probability measure $\mu$ such that for every partition $S_1,S_2,S_3$ of $[d]$, 
        \begin{equation}
            \label{counter_condition}
            \tau_{\mu}(S_1)+\tau_{\mu}(S_2)+\tau_{\mu}(S_3)\le 2m+1.
\end{equation}
\end{theorem}
The following corollary, which follows directly from Theorem \ref{identi_thm}, builds upon the separability condition introduced earlier.
 \begin{corollary}
   \label{sepa_thm}
        Let $\mu$ be defined as in \eqref{eq:MM}. If $N(\mu)\geq 2m-1$, then $\mu$ is identifiable. 
   \end{corollary}
Theorem~\ref{identi_thm} quantifies the contribution of each variable through the diversity index $\textnormal{\textbf{Ind}}_\mu(j)$.
To the best of our knowledge, this is the first result that unifies all previously known identifiability conditions for the model in~\eqref{eq:MM}.
For example, it generalizes the linear independence condition in Theorem~\ref{linearind}, which requires that every variable is $m$-independent and thus guarantees identifiability when $d \ge 3$.
It also extends the result in \cite{vandermeulen_generalized_2022}, which assumes conditional \iid\ variables, while our result only requires conditional independence.
Corollary \ref{sepa_thm} explains why $2m-1$ emerges as a critical threshold for identifiability in existing literature and unifies identifiability conditions from \cite{rabani_learning_2014,tahmasebi_identifiability_2018,vandermeulen_operator_2019}.
Notably, this corollary also resolves a gap in \cite{tahmasebi_identifiability_2018}: whereas their work requires at least $2m$ separable variables to ensure identifiability, our results show that $2m-1$ separable variables suffice. 

Below, we outline the proof of Theorem \ref{identi_thm}; a complete proof is provided in Appendix \ref{app:pf-identify}.  
Our approach is inspired by the Hilbert space embedding technique in \cite{vandermeulen_operator_2019}, which employs a unitary transform connecting the model to the tensor product of Hilbert spaces. Preliminaries on the tensor product of Hilbert spaces are provided in Appendix \ref{tensor_hilbert}.

\noindent\textbf{Proof Sketch.} Consider two probability measures $\mu$ and $\tilde{\mu}$ that are represented in the form of \eqref{eq:MM}. 
Suppose $\mu=\tilde{\mu}$ and $\mu$ satisfies the condition \eqref{condition}. 
There exists a finite measure $\xi$ such that the Radon-Nikodym derivatives $f_{kj} = \frac{d\mu_{kj}}{d\xi},\tilde{f}_{kj}=\frac{d\tilde{\mu}_{kj}}{d\xi}$ are bounded by one, and thus are bounded in $L^2(\xi)$.
Let $f,\tilde{f}$ be the Ranon-Niko derivatives of $\mu,\tilde{\mu}$, respectively.
Applying a unitary transformation (see Lemma \ref{123}), we map $f$ and $\tilde{f}$ to $T$ and $\tilde{T}$, respectively, which reside in the tensor product of Hilbert spaces $L^2(\xi)^{\otimes d}$.
Let $f_{k, S_t} \triangleq \otimes_{j\in S_t}f_{kj}\in L^2(\xi)^{\otimes |S_t|}$ and $\tilde f_{k, S_t} \triangleq \otimes_{j\in S_t}\tilde f_{kj}\in L^2(\xi)^{\otimes |S_t|}$ for $k\in[m]$ and $t=1,2,3$. This allows us to write 
\begin{equation}
\label{eq:T-tildeT}
T=\sum_{k=1}^m (\pi_k f_{k,S_1})\otimes f_{k,S_2}\otimes f_{k,S_3},
\quad
\tilde{T}=\sum_{k=1}^m (\tilde{\pi}_k \tilde{f}_{k,S_1})\otimes \tilde{f}_{k,S_2}\otimes \tilde{f}_{k,S_3}.
\end{equation}
which correspond to the CP decompositions in the tensor product of Hilbert spaces.

Let $f_{S_t}\triangleq (f_{1,S_t},\dots ,f_{m,S_t})\in (L^2(\xi)^{\otimes |S_t|})^{m}$ for $t=1,2,3$. 
Next, we establish a lower bound on the Kruskal rank (see Definition \ref{Kruskal}) of each $f_{S_t}$. 
By Lemma \ref{kru_rank_relation}, the Kruskal rank of $f_{S_t}$ is equal to that of its corresponding Gram matrix $A_{S_t}\in \real^{m\times m}$, where $(A_{S_t})_{kl}=\langle f_{k,S_t},f_{l,S_t}\rangle$. Owing to the inner product structure in Hilbert spaces, the Gram matrix $A_{S_t}$ can be expressed as the \emph{Hadamard product} of the Gram matrices for each variable. 
Specifically, let $f_{j} \triangleq (f_{1j},\dots ,f_{mj})\in (L^2(\xi))^m$ and $A_j$ denote the corresponding Gram matrix. Then,
 \begin{equation}
    \label{connect}
    (A_{S_t})_{kl}  
    =\langle \otimes_{j\in S_t}f_{kj}, \otimes_{j\in S_t}f_{lj}\rangle 
    = \prod_{j\in S_t}\langle f_{kj},f_{lj}\rangle 
    = \prod_{j\in S_t} (A_{j})_{kl}.
\end{equation}
The following crucial lemma demonstrates that the Hadamard product increases the Kruskal rank. 

\begin{lemma}
 \label{Hadamard}
     Suppose $A, B\in\real^{n\times n}$ are real Gram matrices with Kruskal rank $k_A$ and $k_B$ and have no zero main diagonal entries. Then we have
     $$k_{A\circ B}\geq \min\left\{n,k_A+k_B-1\right\}.$$
 \end{lemma}

Prior work \cite[Corollary 5]{horn_rank_2020} establishes a lower bound on the rank of the Hadamard product $A\circ B$, generalizing the classical Schur product theorem \cite[see][Section 7.5]{horn_matrix_2012}.  
In this work, Lemma \ref{Hadamard} extends that result by deriving a lower bound on the \emph{Kruskal rank} of $A\circ B$ tailored to our analysis.
The result also extends the super-additivity property of the Kruskal rank of the Khatri-Rao product, as established in \cite{sid_00}, to general Hilbert spaces.
The proof of Lemma \ref{Hadamard} is provided in Appendix \ref{app:pf-identify}.

By applying Lemma \ref{Hadamard} repeatedly, we deduce that  
\begin{equation}
    \label{eq:Hadamard_rank}
    k_{A_{S_t}}= k_{\circ_{j\in S_t}A_j}\geq \min\left\{m,\sum_{j\in S_t}k_{A_j}-|S_t|+1\right\}.
\end{equation}
Let $k_{f_j}$ denote the Kruskal rank of $f_{j}=(f_{1j},\dots ,f_{mj})$.
By definition, $\sum_{i\in I}a_i f_{ij}\equiv 0$ is equivalent to $\sum_{i\in I}a_i \mu_{ij}\equiv 0$ for every $I\subseteq [m]$.
Therefore, $k_{A_j}=k_{f_j}=\Ind_{\mu}(j)$ and thus $k_{f_{S_t}}=k_{A_{S_t}}\ge \tau_{\mu}(S_t)$. 
Let $f_{k,S_1}'\triangleq (\pi_1 f_{1,S_1},\dots,\pi_k f_{m,S_1})$ with Kruskal rank $k_{f_{S_1}'}$.
Since $\pi_k>0$ for all $k\in[m]$, we have $k_{f_{S_1}'}=k_{f_{S_1}}$.
Combining \eqref{condition} and
\eqref{eq:Hadamard_rank}, 
we obtain that 
\[
k_{f_{S_1}'}+k_{f_{S_2}}+k_{f_{S_3}}
\geq 2m+2.
\]
By applying an extension of Kruskal's theorem in Lemma \ref{Extension} to the tensors in \eqref{eq:T-tildeT}, there exists a permutation $\sigma$ and scalars $C_{1k},C_{2k},C_{3k}$ such that 
\[
\tilde{\pi}_{\sigma(k)} \tilde{f}_{\sigma(k),S_1}=C_{1k} \pi_k f_{k,S_1}, 
\quad \tilde{f}_{\sigma(k),S_2}=C_{2k}f_{k,S_2}, 
\quad \tilde{f}_{\sigma(k),S_3}=C_{3k}f_{k,S_3},
\] 
with $C_{1k}C_{2k}C_{3k}=1$. 
Using the conditions $\int f_{kj} \diff \xi=\int \tilde f_{kj} \diff \xi =1$ and $f_{kj}\ge 0$, we deduce that $C_{2k}=C_{3k}=1$, which implies $C_{1k}=1$. 
Consequently, we conclude that 
\[
f_{kj}=\tilde{f}_{\sigma(k)j},
\quad \pi_k=\tilde{\pi}_{\sigma(k)},
\]
which implies the identifiability result in Theorem \ref{identi_thm}.

Next, we prove the converse result. 
For $d\leq 2m-2$, consider the family of discrete distributions of the form:
\begin{equation}
    \label{Ber_iid}
\mu=\sum_{k=1}^m\pi_k \Bern(\alpha_k)^{\times d}.
\end{equation}
The identifiability of $\mu$ is equivalent to that of binomial mixtures. Specifically, for any $b\in\{0,1\}^{d}$ with $\ell$ nonzero entries, $\mu\{b\}=\sum_{k=1}^m\pi_k\binom{d}{\ell}\alpha_k^\ell(1-\alpha_k)^{d-\ell}$.
Thus, $\mu$ is uniquely determined by $\sum_{k=1}^m \pi_k\alpha_k^j$ for $j\in[d]$, which correspond to the first $d$ moments of the mixing distribution $\sum_{k=1}^m\pi_k \delta_{\alpha_k}$. 
By classical theory of moments, $d\le 2m-2$ moments are insufficient to identify an $m$-atomic distribution  \cite[see, e.g.,][Lemma 30]{wu_optimal_2020}. Hence, $\mu$ is not identifiable. 
Note that $\Ind_{\mu}(j)\leq 2$ for all $j\in[d]$, as any three Bernoulli distributions are linearly dependent.
Consequently, $\tau_{\mu}(S)\le |S|+1$, which implies $\tau_{\mu}(S_1)+\tau_{\mu}(S_2)+\tau_{\mu}(S_3)\leq d+3\leq 2m+1$.

For $d>2m-2$, consider the probablity measure
$
\mu=\sum_{k=1}^m\pi_k \Bern(\alpha_k)^{\times 2m-2}\times \mu_0^{d-2m+2},
$
which reduces the problem to the case $d=2m-2$. 
Here, $\Ind_{\mu}(j)=1$ for $j\ge 2m-1$ and thus $\tau_{\mu}(S)\le |S|+1$ remains valid. 

\section{Rate of Convergence under Incoherence}
\label{sec_roc}
In this section, we focus on the estimation problem of model \eqref{eq:MM}. In the remainder of this paper, we assume each probability measure $\mu_{kj}$ admits a density function $f_{kj}$. The joint density can then be expressed as:
\begin{equation}
    \label{MM_den}
    f(x_1,\dots ,x_d)=\sum_{k=1}^m\pi_k\prod_{j=1}^d f_{kj}(x_j).
\end{equation}
For simplicity, we will henceforth write \eqref{MM_den} as $f=\sum_{k=1}^m\pi_k\prod_{j=1}^d f_{kj}$, with the understanding that the product $\prod_{j=1}^d f_{kj}$ should be interpreted as $\prod_{j=1}^d f_{kj}(x_j)$ unless stated otherwise.

\subsection{Recovering the Component Density: A Perturbation Analysis}
\label{recovery_component_subse}
\par  We say an estimator $\tilde{f}$ is \emph{proper} if it admits the structure \eqref{MM_den}, denoted by $\tilde{f}=\sum_{k=1}^m\tilde{\pi}_k\prod_{j=1}^d\tilde{f}_{kj}$.
We will analyze how the error between  $f$ and $\tilde{f}$ propagates to the components, establishing a perturbation theory that reduces the estimation of model parameters to that of the joint density.  Note that both tasks are harder than the identifiability problem, so we expect stronger conditions than those in Section \ref{sec_identi}.
 We introduce the following incoherence condition in a Hilbert space. 
\begin{definition}[$\mu$-Incoherence]
\label{incoherence}
    Let $f_1,\dots ,f_m$ be elements in a Hilbert space $\mathcal{H}$ and $0\leq \mu< 1$, we say the sequence $\{f_{k}\}_{k=1}^m$ is $\mu$-incoherent if  for any $k\neq k'$,
    $$\left|\langle f_{k},f_{k'}\rangle\right| \leq \mu\|f_k\|_2\|f_k'\|_2.$$
\end{definition}
The above definition has a clear geometric intuition: It can be treated as knowledge of the minimum angle among $f_k$. It is easy to see that $\{f_k\}_{k=1}^m$ is far from parallel as $\mu$ tends to $0$. 
 Based on the incoherence condition, we impose the following technical assumption on the joint density, which is also required for the error analysis of the algorithm proposed later in Section \ref{sec_algo}.
\begin{assumption}[Estimable Condition]
\label{assume}
 For $f=\sum_{k=1}^m\pi_k\prod_{j=1}^{d} f_{kj}$ as in \eqref{MM_den}, we say $f$ is $(\mu,\zeta)$-estimable if 
        \begin{enumerate}
        \item $f_{kj}$'s are square integrable for all $k,j$. For each $j=1,2,\dots ,d$, the set $\{f_{kj}\}_{k=1}^m$ is $\mu$-incoherent with $\mu<1$.
        \item The mixing proportions are uniformly bounded away from zero: $\min_{k\in[m]}\pi_k\geq \zeta>0$.
    \end{enumerate}
\end{assumption}
Now we are ready to present our main result of this subsection, which can be viewed as a robust version of Corollary \ref{sepa_thm}.
\begin{theorem}
\label{Angle_thm}
    Let $f=\sum_{k=1}^m\pi_k\prod_{j=1}^{d} f_{kj}$ be a  $(\mu,\zeta)$-estimable function supported on $[0,1]^d$, and $\tilde{f}=\sum_{k=1}^m\tilde{\pi}_k\prod_{j=1}^d \tilde{f}_{kj}$ be a proper estimator of $f$. Assume that there exists a universal constant $C\geq 1$ such that $\|f_{kj}\|_{\infty},\|\tilde{f}_{kj}\|_{\infty}\leq C$ for all $k,j$. 
    If $\|f-\tilde{f}\|_2 \leq \epsilon$ for $\epsilon<\frac{(1-\mu)^{2m-1}\zeta^2}{32m^{5/2}L_m^2C^{2m}}$, where $L_{m} = 4m^{3/2}(m-1)!>0$, then there exists a permutation $\sigma:[m]\mapsto [m]$,  such that
    \begin{equation*}
        \|f_{kj}-\tilde{f}_{\sigma(k)j}\|_2\leq \frac{8C^2L_{m}}{(1-\mu)^{m-1}\zeta}\epsilon,\quad \|\pi-\sigma(\tilde{\pi})\|_2:=\sqrt{\sum_{k=1}^m(\pi_k-\tilde{\pi}_{\sigma(k)})^2}\leq \frac{16C^{2m-2}L_m^2}{(1-\mu)^{\frac{3(m-1)}{2}}\zeta}\epsilon.
    \end{equation*}
\end{theorem}
Theorem \ref{Angle_thm} shows that under Assumption \ref{assume}, $\|f_{kj}-\tilde{f}_{\sigma(k)j}\|_2$, $\|\pi-\sigma(\tilde{\pi})\|_2$ has the same order as $\|f-\tilde{f}\|_2$.
The result extends the result in \cite{bhaskara_uniqueness_2013,gordon_identification_nodate} to the nonparametric case.
Below, we sketch the proof of Theorem \ref{Angle_thm}. A complete proof is provided in Appendix \ref{Proof_angle}.  

\noindent\textbf{Proof Sketch.}
For $I\subseteq[d]$, let $f_I$ and $\tilde{f}_I$ denote the marginal densities of $f$ and $\tilde{f}$ with respect to the variables indexed by $I$, respectively. Since $f$ and $\tilde{f}$ are supported on $[0,1]^d$, we have $\|f_I-\tilde{f}_I\|_2\leq \|f-\tilde{f}\|_2\leq \epsilon$ from Cauchy-Schwarz inequality. 
In the sequel, we assume without generality that $I=[2m-1]$.  

Similar to the proof of Theorem~\ref{identi_thm}, we represent the joint densities in the tensor product of Hilbert spaces. 
Under the conditions of Theorem \ref{Angle_thm},  $f_{kj},\tilde{f}_{kj}\in L^2([0,1])$ for each $k$ and $j$. 
Thus, by applying a unitary transformation $U$, the joint densities $f$ and $\tilde{f}$ can be represented as finite-rank linear operators $T$ and $\tilde{T}$ in the tensor product space $L^{2}([0,1])^{\otimes (2m-1)}$: 
\begin{equation}
    \label{eq:T_tildeT}
    T=\sum_{k=1}^m\pi_k\otimes_{j=1}^{2m-1}f_{kj},\ 
\quad
\tilde{T}=\sum_{k=1}^m\tilde{\pi}_k\otimes_{j=1}^{2m-1}\tilde{f}_{kj}.
\end{equation} 
Now we consider the mode-1 multiplication of $T$:
For $w\in L^2([0,1])$, we write
$$T\times _1w = \sum_{k=1}^m \pi_k\langle w,f_{k1}\rangle\otimes_{j=2}^{2m-1}f_{kj} \in L^2([0,1])^{\otimes 2m-2}.$$
Then, we unfold $T\times_1 w$ to the following linear operator by a unitary transformation $U'$: 
$$T_{w}=AD_{\pi,w}B^* \in L^2([0,1])^{\otimes (m-1)}\otimes L^2([0,1])^{\otimes (m-1)},$$
where $A=(\otimes_{j=2}^m f_{1j},\dots ,\otimes_{j=2}^m f_{mj}), B=(\otimes_{j=m+1}^{2m-1} f_{1j},\dots ,\otimes_{j=m+1}^{2m-1} f_{mj})$ and $D_{\pi,w}=\diag\{\pi_1$\\
$\langle w,f_{11}\rangle,\dots ,\pi_m\langle w,f_{m1}\rangle\}$. 
Similarly, we map $\tilde{T}$ to $\tilde{T}_w=\tilde{A}D_{\tilde{\pi},w}\tilde{B}^*$. 
Let $\|\cdot\|_{\op}$ denote the operator norm of a linear operator. Since $U$ and $U'$ preserve the inner product, we have $\|T-\tilde{T}\|_{\textnormal{op}}=\|f-\tilde{f}\|_2\leq \epsilon$, $\|T_w-\tilde{T}_w\|_{\op}=\|T\times_1 w-\tilde{T}\times_1w\|_{\op}$.
From the definition of operator norm, we can deduce that
$\sup_{\|w\|_2=1}\|T_w-\tilde{T}_w\|_{\op}\leq \|T-\tilde{T}\|_{\op}\leq \epsilon$. Thus, by Lemma \ref{Weyl},
we obtain the following crucial result:
\begin{equation}
    \label{eq:sing_diff_main}
    \sup_{\|w\|_2=1}\max_{k\in[m]}|\sigma_{k}(T_w')-\sigma_k(\tilde{T}_w')|\leq \epsilon.
\end{equation}

The idea of proof is that if Theorem \ref{Angle_thm} does not hold, we can obtain a lower bound of $|\sigma_{k}(T_w)-\sigma_k(\tilde{T}_w)|$ for some $k\in[m]$ and $w\in L^2([0,1])$ with $\|w\|_2=1$.
Under the incoherence condition, we show that $\sigma_{m}(A),\sigma_{m}(B)\geq \sqrt{\frac{(m-1)!}{(1-\mu)^{m-1}}}$ from Lemma \ref{condition_num_hadamard},
which allows us to focus on the diagonal entries of $D_{\pi,w}$ only.

We first prove that for every $k\in[m]$ and $j\in[d]$, there exists $k'\in [m]$ such that $\|f_{k'j}-\tilde{f}_{kj}\|_2\le \frac{8C^2L_{m}}{(1-\mu)^{m-1}\zeta}\epsilon$. By Lemma \ref{L_12_error} and the assumption on $\epsilon$, it suffices to show that $\sin\theta(f_{k'j},\tilde{f}_{kj})\leq \epsilon'\triangleq \frac{L_m}{(1-\mu)^{m-1}\zeta}\epsilon$.
Suppose on the contrary that there exists some $\tilde{f}_{kj}$ for which $\sin\theta(f_{k'j},\tilde{f}_{kj})>\epsilon'$ for all $k'\in [m]$; without loss of generality, take $j=1$.
Using the probabilistic method, we prove in Lemma \ref{lemma_test} that there exists a test function $w_0\in L^2([0,1])$ with $\|w_0\|_2=1$ such that $|\langle w_0,f_{k'1}\rangle|>\epsilon'\cdot \frac{1}{4m^{3/2}}=\frac{(m-1)!}{(1-\mu)^{m-1}\zeta}\epsilon$ for all $k'\in[m]$, 
yet $\langle w_0,\tilde{f}_{k1}\rangle=0$.
Consequently, $\sigma_{m}(\tilde{T}_{w_0})=0$ whereas $|\sigma_{m}(T_{w_0})|> \sigma_{m}(A)\sigma_{m}(B)\max_{k'\in[m]}|\langle w_0,f_{k'1}\rangle|\geq\epsilon$, which contradicts \eqref{eq:sing_diff_main}.

\par As a result, we build a mapping from $k\in[m]$ to $k'\in[m]$ for each $j\in[d]$, denoted by $\sigma^{(j)}$. Next, we prove the mapping $\sigma^{(j)} :k\mapsto k'$ above  is one-to-one. 
Suppose on the contrary this is not true, then there exists $j\in[d]$ and $k_1,k_2,k'\in[m], k_1\ne k_2$, such that $\|\tilde{f}_{k_1j}-f_{k'j}\|_2,\|\tilde{f}_{k_2 j}-f_{k'j}\|_2\leq 8C^2\epsilon'$; Without loss of generality, take $k'=j=1$.
From the $\mu$-incoherence of $\{f_{k1}\}_{k=1}^m$ and
Lemma \ref{lemma_test}, there exists a test function $w_1\in L^2([0,1])$ with $\|w_1\|_2=1$, such that $|\langle w_1,f_{k1}\rangle|\geq \frac{\sqrt{1-\mu^2}}{4m^{3/2}}$  for all $k\ne 1$, whereas $\langle w_1,f_{11}\rangle=0$.
The latter implies that $|\langle w_1,\tilde{f}_{k_t1}\rangle|=\langle w_1,f_{11}-\tilde{f}_{k_t1}\rangle\leq \|f_{11}-\tilde{f}_{k_t1}\|_2\leq  8C^2\epsilon'$ for $t=1,2$.
Consequently, $|\sigma_{m-1}(T_{w_1})|\geq \frac{\zeta(1-\mu)^{m-1}\sqrt{1-\mu^2}
}{L_m}$, whereas $|\sigma_{m-1}(\tilde{T}_{w_1})|\leq C^{2m-2}\epsilon'$. Combined with the assumption on $\epsilon$, we obtain $|\sigma_{m-1}(T_{w_1})-\sigma_{m-1}(\tilde{T}_{w_1})|>\epsilon$, which contradicts \eqref{eq:sing_diff_main}.

Finally, we prove that $\sigma_j$ are identical for all $j\in[d]$.
Suppose $\sigma_1\ne \sigma_2$.
Then $\sigma_1$ and $\sigma_2$ map two distinct indices $j_1,j_2$ to the same image, say $\sigma_1(1)=\sigma_2(2)=1$.
Define $T'=\sum_{k=1}^m\tilde{\pi}_k f_{\sigma_{1}(k)1}\otimes f_{\sigma_2(k)2}\otimes(\otimes_{j=3}^{2m-1}\tilde{f}_{kj})$.
By the triangle inequality, we deduce that $\|T-T'\|_{\op}\leq 17mC^{2m}\epsilon'$.
Since $\{f_{k1}\}_{k=1}^m$,
$\{f_{k2}\}_{k=1}^m$ are $\mu$-incoherent, 
applying Lemma \ref{lemma_test} again, 
there exist $u,v\in L^2([0,1])$ with $\|u\|_2=\|v\|_2=1$, 
such that $\langle u,f_{11}\rangle=\langle v,f_{12}\rangle=0;|\langle u,f_{k1}\rangle|,|\langle v,f_{k2}\rangle|\geq \frac{\sqrt{1-\mu^2}}{4m^{3/2}}$ for $k=2,3,\cdots m$.
Let $T_{u,v,w}\triangleq T\times_1u\times_2v \times_3 w, T_{u,v,w}'\triangleq T'\times_1u\times_2v\times_3 w$.
Since $\sigma_1(1)=\sigma_2(2)=1$ and $\langle u,f_{11}\rangle=\langle v,f_{12}\rangle=0$, $T_{u,v,w}$ has rank $m-1$, while $T'_{u,v,w}$ has rank at most $m-2$. Treating $T_{u,v,w},T_{u,v,w}'\in L^2([0,1])^{\otimes (2m-4)}$ in the same manner as $T\times_1 w, \tilde{T}\times_1 w$ earlier, we unfold them to $S_{u,v,w},S'_{u,v,w}\in L^2([0,1])^{\otimes (m-2)}\otimes L^2([0,1])^{\otimes (m-2)}$. By choosing $w=\frac{f_{23}}{\|f_{23}\|_2}$, we obtain $|\sigma_{m-1}(S_{u,v,w})-\sigma_{m-1}(S'_{u,v,w})|>17mC^{2m}\epsilon'>\|T-T'\|_{\op}$, which leads to a similar contradiction.

\subsection{Estimation of the Joint Distribution under H\"{o}lder Smoothness Condition}
\label{sec:complexity}

\par In this subsection, our goal is to analyze the complexity of model \eqref{MM_den}. Let $\mathcal{G}_{\mathcal{F}}^{(m,d)}$ be the density class that admits the structure of \eqref{MM_den}, with component densities $f_{kj}$ in class $\mathcal{F}$:
\begin{equation}
\label{MM_class}
    \mathcal{G}_{\mathcal{F}}^{(m,d)}:= \left\{f=\sum_{k=1}^m\pi_k\prod_{j=1}^d f_{kj}:\pi=(\pi_1,...,\pi_m)\in\Delta^{m-1},f_{kj}\in\mathcal{F}\right\}.
\end{equation}
In the following, we will consider a H\"{o}lder smooth density class $\mathcal{F}_{L,q}$ (see Definition \ref{Holder}) for the component densities $f_{kj}$, and derive minimax rate bounds for the class $\mathcal{G}_{\mathcal{F}}^{(m,d)}$ under a suitable metric $\rho$. 
 \begin{theorem}
\label{Holder_rate}
Let $\calF_{L,q}$ denote the class of all $q$-H\"{o}lder smooth densities on $[0,1]$ with smoothness parameter $q$ and constant $L>0$. Given a random sample $X_1,\dots ,X_n\sim f\in\mathcal{G}_{\calF_{L,q}}^{(m,d)}$, we define the minimax risk for class $\mathcal{G}_{\calF_{L,q}}^{(m,d)}$ under a metric $\rho$ as
\begin{equation}
\label{min_max_risk}
        R_{\rho,
    \calF_{L,q}}^*(m,d) \triangleq \inf_{\hat{f}_n}\sup_{f\in\mathcal{G}_{\calF_{L,q}}^{(m,d)}}\Expect[\rho^2(\hat{f}_n,f)].
\end{equation}
Then we have
    \begin{enumerate}
    \item For $n\ge md^{1+\frac{1}{q}}$,
     \begin{equation*}
         \label{mini_He}
         (n\log n)^{-\frac{q}{q+1}}d\lesssim_{L,q} R^*_{H,\calF_{L,q}}(m,d)\lesssim_{L,q} n^{-\frac{q}{q+1}}m^{\frac{q}{q+1}}d.
     \end{equation*}
        \item For all $n\ge 1$,
         \begin{equation*}
             \label{mini_TV}
             (n\log n)^{-\frac{2q}{2q+1}}\lesssim_{L,q} R^*_{\textnormal{TV},\calF_{L,q}}(m,d)\lesssim_{L,q} n^{-\frac{2q}{2q+1}}m^{\frac{2q}{2q+1}}d^{\frac{2q+2}{2q+1}}.
        \end{equation*}  
    \end{enumerate}
  \end{theorem}

We now compare the minimax rates obtained under the latent structure to those for density estimation without latent variables.
It is well known that the minimax rate of estimating a $q$-H\"{o}lder continuous density in $d$ dimensions is of order $n^{-\frac{q}{q+d}}$ in $H$ and $n^{-\frac{q}{2q+d}}$ in $\TV$ \cite[see, e.g.,][Section 32]{Polyanskiy_Wu_2025}, both of which suffer from the curse of dimensionality.
In contrast, Theorem \ref{Holder_rate} shows that the conditional independence structure in our latent variable model retains the minimax behavior of the one-dimensional case, with only a polynomial dependence on $m$ and $d$.
This highlights how leveraging latent structure mitigates the curse of dimensionality in high-dimensional density estimation.
The proof of Theorem \ref{Holder_rate} is based on a classical information-theoretic framework through metric entropy, and the detail is provided in Appendix \ref{pf:stat_limit}.

\section{Algorithm for Recovery of the Components}
\label{sec_algo}
\subsection{An Operational Method for Recovery}
In this subsection, we will develop an operational procedure for recovering each component density $f_{kj}$ from an estimator of the joint density $f$ in model \eqref{MM_den}. We propose a recovery algorithm based on the simultaneous diagonalization method introduced by \cite{leurgans_decomposition_1993}. This method has been applied in some special cases of model \eqref{MM_den} in earlier works. \cite{bonhomme_estimating_2016} applied the technique to density estimation by projecting the component densities onto the top terms of an (infinite) orthogonal basis and estimating their coefficients from a random sample. \cite{gordon_identification_nodate} applied the same method to the Bernoulli mixture model and analyzed the robustness of the algorithm.
\par
We focus on the case that the joint density $f$ satisfies Assumption \ref{assume}. 
We first consider the case $d=2m-1$, the smallest dimension that ensures identifiability.  
We present the recovery procedure in Algorithm 1 below.
A more detailed discussion of Algorithm 1 is provided in Appendix \ref{algo_compare}.
\begin{algorithm}
\label{recovering_algo}
    \caption{Recover the component density from the estimator of joint density}
    \begin{algorithmic}[1]
    \Require An estimator $\hat{f}$ for the density $f=\sum_{k=1}^m \pi_k\prod_{j=1}^{2m-1} f_{kj}$ on $[0,1]^{2m-1}$ 
    \Ensure $\hat{f}_{k1}$ for $k=1,2,\dots ,m$
    \State \label{line_1}  Calculate $\hat{T}_{+}(y,z) = \int \hat{f}(y,z,x_{2m-1}) dx_{2m-1}$, where $y=(x_1,\dots ,x_{m-1})$ and $z=(x_m,\dots ,x_{2m-2})$.    
    \State \label{line_2} Let $\hat{T}_{+,m}(y,z) = \textnormal{argmin}_{\rank(T) \le m}\|T - \hat{T}_{+}\|_{\textnormal{op}} = \sum_{k=1}^m \hat{\lambda}_k\hat{\phi}_{k}(y)\hat{\psi}_{k}(z)$, the top $m$ truncation of singular value decomposition (SVD) 
     \State \label{line_3} Choose some subset $A\subset [0,1]$
   \For{$l,t = 1, 2, \dots, m$} \State\label{line_4}
    $\hat{\eta}_{lt} \gets \frac{1}{\hat{\lambda}_t} \int_A \hat{\phi}_l(y) \hat{f}(y,z,x_{2m-1}) \hat{\psi}_t(z)  \textnormal{d}y  \textnormal{d}z  \textnormal{d}x_{2m-1}$
    \EndFor
     \State \label{line_5} Let $\hat{\eta}_A=(\hat{\eta}_{lt})_{m\times m}$, calculate $\hat{W}\gets (\hat{w}_1,\dots ,\hat{w}_m)$ where $\hat{w}_1,\dots ,\hat{w}_m$ are $L_2$ unit eigenvectors of $\hat{\eta}_A$
  \For{$k=1,2,\dots ,m$}   \State\label{algo_return}
     $\hat{g}_{k}(y)\gets \sum_{h=1}^m \hat{w}_{kh}\phi_h(y)$, $\hat{h}_{k}\gets \hat{g}_k/\|\hat{g}_k\|_1$, $\hat{f}_{k1}\gets\int \hat{h}_k \textnormal{d}x_2\dots \textnormal{d}x_{m-1}$ 
     \EndFor
    \end{algorithmic}
\end{algorithm} 

Now we show that Algorithm 1 correctly recovers the component density under Assumption \ref{assume} given a good choice of subset $A$.

\begin{theorem}[Correctness of Algorithm 1]
\label{algo_error}
   Suppose the density function $f=\sum_{k=1}^m \pi_k\prod_{j=1}^{2m-1}f_{kj}$ on $[0,1]^{2m-1}$ is $(\mu,\zeta)$-estimable, and $\|f_{kj}\|_{\infty}\leq C$ for all $k,j$. Suppose the following conditions hold:
    \begin{enumerate}
        \item The Lebesgue measure of $A$ is large: $\mu_{Leb}(A)\geq \mu_0$.
        \item $a_k = \int_Af_{k(2m-1)}(x)dx$ are lower bounded and well separated:
        $$\min_{k\in [m]}a_k\geq\delta,\min_{k\ne k'}|a_k-a_k'| \geq\delta.$$
    \end{enumerate}
   Then for a density estimator $\hat{f}$ satisfying $\|\hat{f}-f\|_2\leq \epsilon$ for some $\epsilon< \frac{\zeta(1-\mu)^m}{4(m-1)!}$, Algorithm 1 outputs $\hat{f}_{k1}$ such that
    \begin{equation*}
        \|\hat{f}_{k1}-f_{\sigma(k)1}\|_2 \leq \frac{L_{C,m}\epsilon}{\zeta^3(1-\mu)^{3m} \delta\sqrt{\mu_{0}}}
    \end{equation*}
    for a permutation $\sigma: [m]\mapsto [m]$ and a universal constant $L_{C,m}>0$ depending on $C$ and $m$ only.
\end{theorem}
\begin{remark}
If each $f_{kj}$ is a probability mass function supported on the discrete set $\{1, 2, \dots, N\}$, then Algorithm 1 can still be applied with minor modifications. Specifically, the integrals in Algorithm 1 should be replaced with summations, and the random set $A$ should be sampled as a random weight vector over ${1, 2, \dots, N}$. According to prior results in \cite{bhaskara_smoothed_2014}, under the incoherence condition, Condition 2 in Theorem \ref{algo_error} is satisfied with probability 1, and the parameter $\delta$ will depend on the incoherence level $\mu$. In this discrete setting, the error bound will incur an additional factor that depends only on $N$.
\end{remark}
Theorem \ref{algo_error} establishes that, as long as $\hat{f}$ is sufficiently close to $f$, we can accurately recover each component density $f_{k1}$ for $k = 1, 2, \dots, m$. Notably, the theorem relies only on the incoherence condition, rather than the stronger linear independence condition often assumed in previous work. In the general case where $d \geq 2 m-1$, we can repeatedly apply our algorithm to submodels of size $2 m-1$ to recover all component densities $f_{kj}$ for every $k$ and $j$, requiring $d$ such repetitions. The proof of Theorem \ref{algo_error} is provided in Appendix \ref{proof_algo}.
\subsection{Simulations}
\label{simulations}
We set up two simulations for the case where $f_{kj}$'s are probability mass functions. The first simulation is the conditional \iid\  model in Example~\ref{ex:CI}, and the second is for the Bernoulli mixture model in Example~\ref{ex:BMM}. In both simulations, we set $m=3,d=5$, so the true probability mass is $f=\sum_{k=1}^3 \pi_k \prod_{j=1}^5 f_{kj}$.
We report the following measure $e = \sum_{k=1}^m \|f_{k1}-\hat{f}_{k1}\|_2.$ To obtain $\hat{f}$, we will first draw a random sample $X_1,\dots ,X_n\sim f$, and use empirical estimate. To control the error between $\hat{f}$ and $f$, we set an exponential growth for sample size $n=2^{17},\dots ,2^{24}$. The experiment is repeated 10 times, and we report the mean and variance of error $e$ by a log-log plot.

\begin{figure}[ht]
    \centering
    \subfloat[Error plot for conditional i.i.d. model]{%
        \includegraphics[width=0.48\linewidth]{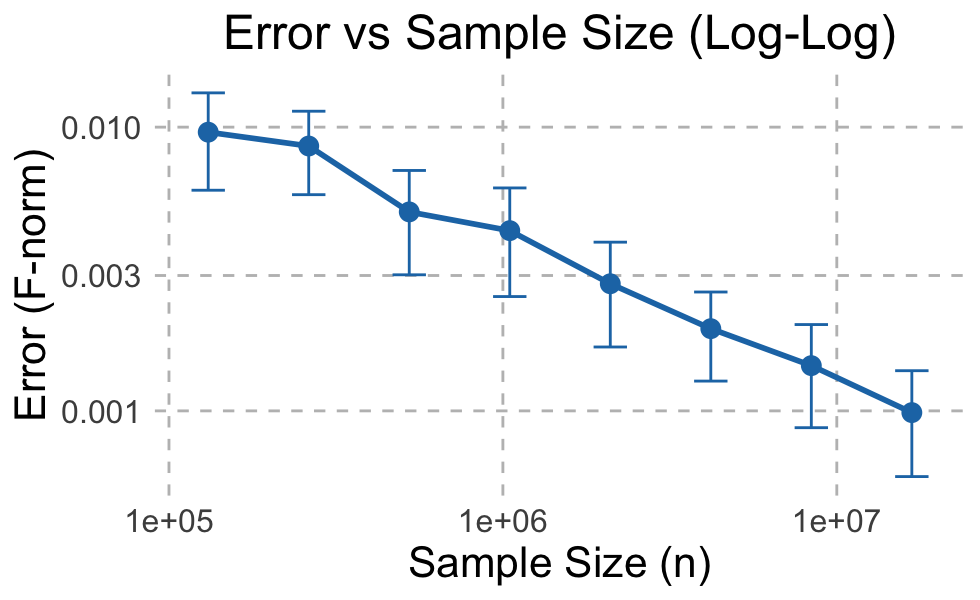}
        \label{CI_plot}
    }
    \hfill
    \subfloat[Error plot for Bernoulli mixture model ]{%
        \includegraphics[width=0.48\linewidth]{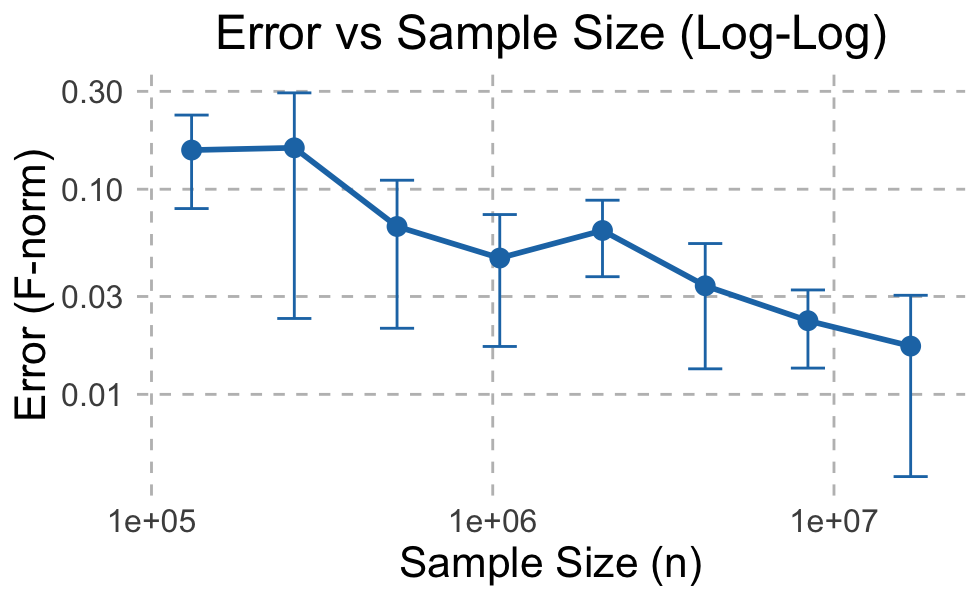}
        \label{BMM_plot}
    }
\end{figure}

\textit{Simulation study 1: Conditional i.i.d. model.} We set the support of $f_{kj}$'s as $\{1,2,3,4\}$, and the probability mass function can be represented by a $4$-dim vector. We set $f_{1}=f_{11}=\dots =f_{15} = (\frac{1}{4},\frac{1}{4},\frac{1}{4},\frac{1}{4})$; $f_{2}=f_{21}=\dots =f_{25} = (0,0,\frac{1}{2},\frac{1}{2})$ ; $f_{3}=f_{31}=\dots =f_{35} = (\frac{1}{2},\frac{1}{2},0,0)$. The mixing proportion $\pi=(0.2,0.3,0.5)$.  The result is shown in Figure \ref{CI_plot}.

\textit{Simulation study 2: Bernoulli mixture model }. For $f_{kj}\sim \Bern(\alpha_{kj})$, we set $\alpha_{kj} = 0.1j+0.2(k-1)$ and the mixing proportion to be $\pi=(0.2,0.3,0.5)$.  The result is shown in Figure \ref{BMM_plot}.
\par Now we discuss the simulation results. First, as the sample size increases, the log error of the component density exhibits a clear linear decay. Since the error of $\hat{f}$ and $f$ has rate $n^{-c}$ with high probability, this experiment confirms the linear relationship between the joint density error and the component density error, as stated in Theorem \ref{algo_error}. Notably, in both simulations, the linear independence condition is not required. The superior performance of the conditional i.i.d. model compared to the Bernoulli mixture model can be attributed to its lower number of parameters and a better separation of the true parameters.

\section{Discussion}
This paper proposes a high-dimensional nonparametric latent structure model. We introduce an identifiability theorem that unifies existing
conditions. In particular, we demonstrate that the increasing dimensionality, coupled with diversity in variables, is beneficial to the identifiability. We also establish a perturbation theory under incoherence and derive minimax risk bounds for high-dimensional nonparametric density estimation, which add up to quantitative rates of convergence. We also develop a recovery algorithm from an estimator of the joint density, which can successfully recover the component densities under incoherence. 
\par There are also some problems to be further investigated under our model:
\begin{itemize}
    \item \textit{Identifiability conditions.} For now, Theorem \ref{identi_thm} is built on a $3$-partition of $[d]$. Such a condition could be replaced by properties only depending on $\mu$. Besides, the condition is still not necessary.
    \item \textit{Full use of diversity.} For large $d$, we estimate the component only using $2m-1$ variables. Using more variables could be more beneficial.
\end{itemize}

\section*{Acknowledgment}
The authors thank Anru Zhang for helpful discussions at the onset of the project. 
The authors are also grateful to anonymous reviewers for helpful comments.

\bibliography{COLt2025.bib}
\bibliographystyle{alpha}
\newpage
\appendix


\section{Proof in Section \ref{sec_identi}}
\label{proof_identifiability}
\subsection{Tensor of Hilbert spaces}
\label{tensor_hilbert}
We first establish the framework of the tensor of Hilbert spaces. Here, we only introduce the definitions and propositions we need to avoid the ambiguity of the notations we use. Proofs of classical results are omitted in this subsection; see Chapter 2 of \cite{reed1980methods,kadison1983fundamentals} for details.
\par Let $\mathcal{H}$,$\mathcal{H}'$ be two Hilbert spaces with basis $\{e_n\}_{n=1}^{\infty},\{e_n'\}_{n=1}^{\infty}$  and inner product $\langle\cdot,\cdot\rangle_\mathcal{H},\langle \cdot,\cdot\rangle_{\mathcal{H'}}$. For $h\in\mathcal{H},h'\in\mathcal{H} ' $, let $h\otimes h'$ (also called a simple tensor) be the bilinear form acting on $\mathcal{H}\times \mathcal{H} ' $: For $g\in\mathcal{H}, g'\in\mathcal{H}'$,
\begin{equation}
    \label{bilinear}
    h\otimes h' (g,g') := \langle h,g\rangle_\mathcal{H} \langle h',g'\rangle_{\mathcal{H'}}.
\end{equation}
Let $\mathcal{E}=\Span \{h\otimes h': h\in\mathcal{H},h'\in\mathcal{H}'\}$ be the linear combinations of all bilinear forms. The tensor of Hilbert spaces $\mathcal{H}$ and $\mathcal{H}'$, denoted by $\mathcal{H}\otimes \mathcal{H}'$, is defined by the completion of $\mathcal{E}$. It can be verified that (See e.g., Proposition 2 in Chapter 2 of \cite{reed1980methods})
$\mathcal{H}\otimes \mathcal{H}'$ is a Hilbert space with basis $\{e_n\otimes e_m'\}_{n,m=1}^{\infty}$ and the following inner product rule:
$$\langle e_i\otimes e'_{j}, e_k\otimes e'_l\rangle_{\mathcal{H}\otimes \mathcal{H}'} = \delta_{ik}\delta_{jl}.$$
Under this rule, it can be verified that the inner product of two simple tensors is
\begin{equation}
    \label{inner}
    \langle h_1\otimes h_2,h_1'\otimes h_2'\rangle =\langle h_1,h_2\rangle_{\mathcal{H}} \langle h_1', h_2'\rangle_{\mathcal{H}'}.
\end{equation}
Note that the definition of inner product from equation \eqref{inner} is equivalent to the one defined on the basis, so we will use \eqref{inner} later on. Now we turn to the tensor product of $d$ Hilbert spaces $\mathcal{H}_1\otimes\dots \otimes \mathcal{H}_d$. By Proposition 2.6.5 in \cite{kadison1983fundamentals}, we know that the tensor product is associative in the sense of isomorphism. Thus, $\mathcal{H}_1\otimes\dots \otimes \mathcal{H}_d$ is defined as the completion of the span of order-$d$ simple tensors $span \{h_1\otimes\dots \otimes h_d: h_i\in H_i,\quad i=1,2,\dots ,d\}$, with the inner product 
$$\langle h_1\otimes\dots h_d,h_1'\otimes\dots \otimes h'_d\rangle_{\mathcal{H}_1\otimes\dots \otimes\mathcal{H}_d} = \langle h_1,h_1'\rangle_{\mathcal{H}_1}\dots \langle h_d,h_d'\rangle_{\mathcal{H}_d}.$$
For a Hilbert space $\mathcal{H}$, the notation $\mathcal{H}^{\otimes d}$ is defined as the $d$-tensor power of $\mathcal{H}$, i.e., $\mathcal{H}^{\otimes d} = \underbrace{\mathcal{H}\otimes\dots \otimes \mathcal{H}}_{\text{$d$ times}}$. In the remainder, the notations should be viewed as the definitions above.
\par Tensor of Hilbert spaces $\mathcal{H}_1\otimes\dots \otimes\mathcal{H}_d$ has a natural isomorphism to the product of Hilbert spaces $\mathcal{H}_1\times\dots \times \mathcal{H}_d$, like the unfolding of a high order tensor in the Euclidean space.  The following classical result reveals the relationship in $L^2$ space (See e.g., Theorem II.10 (a) in \cite{reed1980methods}, also Lemma 5.2 in \cite{vandermeulen_operator_2019}).
\begin{lemma}
 \label{123}
    For a measurable space $(\Psi,\mathcal{G},\gamma)$, there exists a unitary transform $U: L^2(\Psi,\mathcal{G},\gamma)^{\otimes d}\to L^2(\Psi^{\times d},\mathcal{G}^{\times d},\gamma^{\times d})$ such that for all $f_1,\dots ,f_d\in L^2(\Psi,\mathcal{G},\gamma)$,
    \begin{equation}
        \label{unitary}
        U(f_1\otimes\dots \otimes f_d) = f_1(\cdot)\dots f_d(\cdot).
    \end{equation}
\end{lemma}

\subsection{Proof of Theorem \ref{identi_thm}}
\label{app:pf-identify}
Before proving Theorem \ref{identi_thm}, we need to formally define the Kruskal rank:
\begin{definition}[Kruskal rank of a matrix]
\label{Krank_mat}
    Let $M\in\real^{m\times n}$ be a real matrix. The Kruskal rank of $M$ is defined as the maximum number $k$ such that any $k$ columns of $M$ are linearly independent. Denote the Kruskal rank of $M$ by $k_M$.
\end{definition}
\begin{definition}[Kruskal rank in Hilbert spaces]
   \label{Kruskal}
       Let $h=(h_1,\dots,h_m)\in\mathcal{H}^m$. We say $h$ is $k$-independent if, for any size-$k$ index set $S=\{i_1,\dots,i_k\}\subseteq [m]$, $h_{i_1},\dots,h_{i_k}$ are linearly independent. The Kruskal rank of $h$ is the maximum number $k$ such that $h$ is $k$-independent. Denote the Kruskal rank of $h$ by $k_h$.
  \end{definition}
The following lemma reduces the analysis of general Hilbert spaces to the associated Gram matrices. 
 \begin{lemma}
 \label{kru_rank_relation}
      Let $h=(h_1,\dots,h_n)\in\mathcal{H}^n$, and let $G=(\langle h_i,h_j\rangle)_{i,j=1}^n\in\real^{n\times n}$ denote the associated Gram matrix. Then, the Kruskal ranks satisfy $k_h=k_G$. 
 \end{lemma}
 \begin{proof}
      We first prove $k_G\leq k_h$. By the definition of Kruskal rank, there exist $k_h+1$ elements in $h$ that are linearly dependent. Without loss of generality, assume these are $h_1,\dots,h_{k_h+1}$. 
      Partition the Gram matrix $G$ into blocks: 
    $$G=\begin{pmatrix}
        G_{11}& G_{12}\\
        G_{21}& G_{22}
    \end{pmatrix},$$
  where $G_{11}\in\real^{(k_h+1)\times (k_h+1)}$ is the submatrix corresponding to the inner products of $h_1,\dots,h_{k_h+1}$. 
  Since these elements are linearly dependent, $G_{11}$ is rank deficient. 
  By the row inclusion property \cite[see][Observation 7.1.12]{horn_matrix_2012}, the first $k_h+1$ columns of $G$ are linearly dependent. Thus, $k_G\leq k_h$.

Next, we prove $k_G\geq k_h$.
By the definition of Kruskal rank, every subset of $k_h$ elements in $\{h_1,\dots,h_n\}$ is linearly independent.
Consequently, every principal submatrix of $G$ of order $k_h$ has full rank. 
Applying the row inclusion property again, any $k_h$ columns of $G$ are linearly independent. 
Therefore, $k_G\geq k_h$.
\end{proof}

\begin{proof}[Proof of Lemma \ref{Hadamard}]
We prove two cases separately. 

\textbf{Case 1}: $k_A+k_B\ge n+1$. We prove $A\circ B$ is positive definite, which implies $k_{A\circ B}\ge n$. 
Suppose $x^{\top}(A\circ B)x=0$. 
Using the factorization $A=P^{\top}P,B=Q^{\top}Q$, where $P=A^{1/2}$, $Q=B^{1/2}$, we compute:
\begin{align*}
    0=x^{\top}(A\circ B)x  &=  \tr\left(A\diag(x)B\diag(x)\right)\\
    &= \tr\left(P^{\top}P\diag(x)Q^{\top}Q\diag(x)\right)\\
    &=\tr(P\diag(x)Q^{\top}Q\diag(x)P^{\top}) 
    = \| P\diag(x)Q^{\top}\|_F^2.
\end{align*}
This implies $P\diag(x)Q^{\top}=0$.
Let $P=(p_1,\dots,p_n),Q=(q_1,\dots,q_n)$, where $p_i,q_i$ are columns vectors.
Then
\[
C=P\diag(x)Q^{\top}=\sum_{i=1}^{n}x_ip_iq_i^{\top}.
\]
Since $A,B$ no zero diagonal entries, $p_i\ne 0$ and $q_i\ne 0$ for all $i\in [n]$.

By Lemma \ref{kru_rank_relation}, $k_Q=k_B$, so $q_1,\dots,q_{k_B}$ are linearly independent.
For each $j=1,\dots,k_B$, let $\calV_j=\Span\{q_1,\dots,q_{j-1},q_{j+1},\dots,q_{k_B}\}$ and project $q_j$ onto the orthogonal complement $\calV_j^\perp$ denoted by $\Pi_{\calV_j^\perp}(q_j)$. 
By linear independence, $q_j\not\in \calV_j$ and thus $w_j\triangleq \Pi_{\calV_j^\perp}(q_j) \ne 0$.
By construction, $q_j^\top w_j\ne 0$ and $q_i^\top w_j=0$ for $i\ne j \le k_B$. 
Therefore,
\[
0=C w_j =(x_j q_j^\top w_j) p_j +\sum_{i=k_B+1}^{n} (x_i q_i^\top w_1) p_i.
\]
Since $k_A\ge n-k_B+1$ and $k_P=k_A$ by Lemma \ref{kru_rank_relation}, the vectors $p_j,p_{k_B+1},\dots,p_n$ are linearly independent. 
Then, $x_j q_j^\top w_j=0$ and thus $x_j=0$.

Since $q_i\ne 0$ for $i\in[n]$, the union of hyperplanes $\cup_{i=1}^n \{w:q_i^\top w=0\}$ has Lebesgue measure zero.
Hence, there exists $w\in\real^n$ such that $q_i^\top w\ne 0$ for all $i\in [n]$. 
Therefore,
$$0=C w = \sum_{i=k_B+1}^{n} (x_i q_i^\top w) p_i.$$
Since $p_{k_B+1},\dots,p_n$ are linearly independent, it follows that $x_i q_i^\top w=0$ and thus $x_i=0$ for $i=k_B+1,\dots, n$.
We obtain $x=0$ and conclude that $A\circ B$ is positive definite. 

\textbf{Case 2}: $k_A+k_B\leq n$. 
We prove that every principal submatrix of $A\circ B$ of order $m\triangleq k_A+k_B-1$ is nonsingular. 
By the row inclusion property of positive semi-definite matrices \cite[see][Observation 7.1.12]{horn_matrix_2012}, this implies every $m$ columns of $A\circ B$ are linearly independent.
Let $C'=A'\circ B'$ denote an arbitrary principal submatrix of $A\circ B$ of order $m$. 
Since $k_A,k_B\ge 1$ due to the nonzero diagonals, we have $m=k_A+k_B-1\ge \max\{k_A,k_B\}$. 
The Kruskal ranks are inherited by those principal submatrices:
\begin{align*}
    &\text{every } k_A \text{ columns of } A \text{ are linearly independent }\\
    \implies & \text{every principal submatrix of } A \text{ of order $k_A$ has full rank }\\
    \implies & \text{every principal submatrix of } A' \text{ of order $k_A$ has full rank }\\
    \implies & \text{every } k_A \text{ columns of }  A' \text{ are linearly independent}. 
\end{align*}
It follows that $k_{A'}\geq k_A$. Similarly, $k_{B'}\geq k_B$. 
The submatrices $A'$ and $B'$ are positive semidefinite with no zero diagonals, and their Kruskal ranks satisfy $k_{A'}+k_{B'}\ge k_A+k_B=m+1$. Since $C'$ is a matrix of order $m$, Case 1 implies that $C'$ has full rank.
\end{proof}

The following lemma \cite[Theorem 5.1]{vandermeulen_generalized_2022} is an adaptation of Kruskal's theorem in the tensor of Hilbert spaces.
\begin{lemma}[Hilbert space extension of Kruskal's theorem]
\label{Extension}
    Let $x=(x_1,\dots,x_m)\in \mathcal{H}_1^m$, $y=(y_1,\dots,y_m)\in \mathcal{H}_2^m$, and $ z=(z_1,\dots,z_m)\in \mathcal{H}_3^m$ have Kruskal ranks $k_x,k_y$ and $k_z$, respectively.
    Suppose that $k_x+k_y+k_z\geq 2m+2$.
    If $a=(a_1,\dots,a_m)\in \mathcal{H}_1^m$, $b=(b_1,\dots,b_m)\in \mathcal{H}_2^m$, $c=(c_1,\dots,c_m)\in \mathcal{H}_3^m$, and
    \begin{equation*}
        \sum_{k=1}^m x_k\otimes y_k\otimes z_k = \sum_{k=1}^m a_k\otimes b_k\otimes c_k,
    \end{equation*}
    then there exists a permutation $\sigma:[m]\to [m]$ and $D_x,D_y,D_y\in \real^m$ s.t. $a_{\sigma(k)}=x_kD_{x}(k),b_{\sigma(k)}=y_kD_{y}(k)$ and $c_{\sigma(k)}=z_kD_{z}(k)$ with $D_x(k)D_y(k)D_z(k)=1$ for all $k\in[m]$.
\end{lemma}
Now we are ready to prove Theorem \ref{identi_thm}. 
\begin{proof}[Proof of Theorem \ref{identi_thm}]
    For two joint probability measure $\mu,\tilde{\mu}$ having the form as model \eqref{eq:MM}, suppose $\mu=\tilde{\mu}$ with parameters $(\pi_k,\mu_k),(\tilde{\pi}_k,\tilde{\mu}_k)$, and $\mu$ satisfies the condition in the statement of Theorem \ref{identi_thm}. Define the finite measure
    \begin{equation*}
    \xi = \sum_{k,j}(\mu_{kj}+\tilde{\mu}_{kj}).
    \end{equation*}
    Then the Radon-Nikodym derivatives $f_{kj} = \frac{\diff\mu_{kj}}{\diff\xi},\tilde{f}_{kj}=\frac{\diff\tilde{\mu}_{kj}}{\diff\xi}$ are bounded by $1$, thus $f_{kj},\tilde{f}_{kj} \in L^1(\real,\mathcal{B}(\real),\xi)\cap L^2(\real,\mathcal{B}(\real),\xi)$ for all $k,j$. As a consequence, the density functions of $\mu$ and $\tilde{\mu}$ with respect to $\xi^{\times d}$ have the form
\begin{equation*}
    \label{density}
    f(x_1,\dots ,x_d)=\sum_{k=1}^m \pi_k\prod_{j=1}^d f_{kj}(x_j),\quad\tilde{f}(x_1,\dots ,x_d)=\sum_{k=1}^m \tilde{\pi}_k\prod_{j=1}^d \tilde{f}_{kj}(x_j).
\end{equation*}
For simplicity, we will write $f_{kj}(x_j)$ as $f_{kj}$ if the notation has no ambiguity. We now rearrange $f$ and $\tilde{f}$ along the partition $S_1,S_2,S_3$ of $[d]$:
$$f=\sum_{k=1}^m\pi_k \prod_{i\in S_1}f_{ki}\prod_{j \in S_2}f_{kj}\prod_{l\in S_3}f_{kl},\quad\tilde{f}=\sum_{k=1}^m\tilde{\pi}_k \prod_{i\in S_1}\tilde{f}_{ki}\prod_{j \in S_2}\tilde{f}_{kj}\prod_{l\in S_3}\tilde{f}_{kl}.$$
 Now, applying Lemma \ref{123}, there exists a unitary transform $U:L^2(\real,\mathcal{B}(\real),\xi)^{\otimes d}\to L^2(\real^{d},\mathcal{B}(\real)^{d},\xi^{\times d})$ such that \eqref{unitary} holds. Now, by linearity of $U^{-1}$ we have
\begin{equation*}
    \label{tensorA}
    T=U^{-1}(f) = \sum_{k=1}^m(\pi_k \otimes_{i\in S_1}f_{ki})\otimes(\otimes_{j \in S_2}f_{kj})\otimes (\otimes_{l\in S_3}f_{kl}),
\end{equation*}
and
\begin{equation*}
    \tilde{T}=U^{-1}(\tilde{f})=\sum_{k=1}^m
    (\tilde{\pi}_k \otimes_{i\in S_1}\tilde{f}_{ki})\otimes(\otimes_{j \in S_2}\tilde{f}_{kj})\otimes(\otimes_{l\in S_3}\tilde{f}_{kl}).
\end{equation*}
From $\mu=\tilde{\mu}$ we know $f=\tilde{f}$, thus $T=\tilde{T}$. We only need to show $f_{kj}=\tilde{f}_{kj}$ up to a permutation from $T=\tilde{T}$. Let $f_{k,S_t}:=\otimes_{i\in S_t}f_{ki}$ for simplicity and $f_{S_t}=(f_{1,S_t},\dots ,f_{m,S_t})$ for $t=1,2,3$. Similarly, we define $\tilde{f}_{k,S_t}$ and $\tilde{f}_{S_t}$ for $\tilde{f}$. From Lemma \ref{kru_rank_relation} and Lemma \ref{Hadamard}, we have the following lower bound for the Kruskal rank of $f_{S_1}$:
\begin{align*}
    k_{f_{S_1}} &= k_{A_{S_1}} = k_{\circ_{j\in S_1}A_j}\\
    &\geq \min\{m,\sum_{j\in S_1}k_{A_j}-|S_1|+1\}\\
    &\geq\min \{m, \sum_{j\in S_1}\textnormal{\textbf{Ind}}_\mu(j)-|S_1|+1\} = \min\{m,\textnormal{\textbf{Ind}}_\mu(S_1)-|S_1|+1\} = \tau_{\mu}(S_1),
\end{align*}
where $A_{S_1},A_j$ is defined as in \eqref{connect}. Similarly, for $k_{f_{S_2}}$ and $k_{f_{S_3}}$ we have $k_{f_{S_2}}\geq \tau_{\mu}(S_2),k_{f_{S_3}}\geq \tau_{\mu}(S_3)$. Now from the condition \eqref{condition}, applying Lemma \ref{Extension} for $A$ and $B$, we conclude that there exists a permutation $\sigma: [m]\to[m]$ and $D_{S_1},D_{S_2},D_{S_3}\in\real^m$, such that for all $k\in[m]$, $D_{S_1}(k)D_{S_2}(k)D_{S_3}(k)=1$ and
$$ \tilde{\pi}_{\sigma(k)}\tilde{f}_{\sigma(k),S_1} = \pi_kD_{S_1}(k)f_{k,S_{1}},\ \tilde{f}_{\sigma(k),S_t} = D_{S_t}(k)f_{k,S_{t}}, t=2,3.$$
Applying the unitary transform $U$ on them, we have
$$\tilde{\pi}_{\sigma(k)}\prod_{j\in S_1}\tilde{f}_{\sigma(k)j}=\pi_kD_{S_1}(k)\prod_{j\in S_1}f_{kj}, \prod_{j\in S_t}\tilde{f}_{\sigma(k)j} = D_{S_t}(k)\prod_{j\in S_t}f_{kj}.$$
Since $f_{kj},\tilde{f}_{kj}$ are all density functions, we know $D_{S_t}(k)=1$ for all $k$ and $t=2,3$. Thus, from $D_{S_1}(k)D_{S_2}(k)D_{S_3}(k)=1$ we know $D_{S_1}(k)=1$ for all $k$ as well, which implies $\pi_k=\tilde{\pi}_{\sigma(k)}, f_{kj}=\tilde{f}_{\sigma(k)j}$ for all $k,j$. Now for any measurable set $A\in \Psi$, $\mu_{kj}(A) = \int_A f_{kj}d\xi = \int_A \tilde{f}_{\sigma(k)j}d\xi =\tilde{\mu}_{\sigma(k)j}(A)$, which implies $\mu_{kj}=\tilde{\mu}_{kj}$, as desired.  
\par Now it remains to find a $\mu_0$ such that \eqref{counter_condition} holds but not identifiable. Here we consider two mixtures of binomial distribution $\mu_0=\sum_{k=1}^m\pi_k\mu_{k}^{\times 2m-2}$ and $\tilde{\mu}_0=\sum_{k=1}^m\tilde{\pi}_k\tilde{\mu}_{k}^{\times 2m-1}$ with $d=2m-2$, where $\mu_{k}\sim \Bern(\alpha_k),\tilde{\mu}_{k}\sim \Bern(\beta_k)$. We will construct $\mu_0,\tilde{\mu}_0$, such that $\mu_0$ satisfies condition \eqref{counter_condition}, $\mu_0=\tilde{\mu}_0$, but $\mu_{k}\neq\tilde{\mu}_k$ by a permutation.
\par Let $\pi_k = \frac{1}{2^{2m-1}}\binom{2m-1}{2k-2}$ and $\tilde{\pi}_{k}=\frac{1}{2^{2m-1}}\binom{2m-1}{2k-1}$ for $k=1,2,\dots,m$. Then $\sum_{k=1}^m\pi_k=\sum_{k=1}^m\tilde{\pi}_k=1$. For all $k\in  [m]$, let $\alpha_k= c(2k-2), \beta_k=c(2k-1)$, where $c>0$ is a small constant s.t. $\alpha_k,\beta_k\in [0,1]$. 

We first show $\mu_0$ satisfies \eqref{counter_condition}. From $\alpha_k\ne \alpha_k'$ for $k\ne k'$, we know that $\{\mu_k\}_{k=1}^m$ is $2$-independent but not $3$-independent. Thus, for $m\geq 3$ and the partition $S_1=\{1,\dots,m-2\},S_2=\{m-1,\dots,2m-3\},S_3=\{2m-2\}$ of $[2m-2]$, we have 
$$\sum_{t=1}^3 \tau_{\mu_0}(S_t) = \sum_{t=1}^3\min\{m,\sum_{j\in S_t}\textnormal{\textbf{Ind}}(j)-|S_t|+1\}= \sum_{t=1}^3\min\{m,|S_t|+1\} = 2m+1.$$
Now we show that $\mu_0=\tilde{\mu}_0$ to complete the proof. For any $a=(a_1,...,a_{2m-2})\in\{0,1\}^{2m-2}$, suppose $\|a\|_0:=\#\{i:a_i\neq 0\}=l\leq 2m-2$, we have
\begin{align*}
    \mu_{0}(a)-\tilde{\mu}_0(a) &= \frac{1}{2^{2m-1}}\sum_{k=1}^m\left(\binom{2m-1}{2k-2}\alpha_k^{l}(1-\alpha_k)^{2m-2-l} - \binom{2m-1}{2k-1}\beta_k^{l}(1-\beta_k)^{2m-2-l}\right)\\
    &= \frac{1}{2^{2m-1}}\sum_{k=1}^m\sum_{s=0}^{2m-2-l}(-1)^s\left(\binom{2m-1}{2k-2}\alpha_k^s-\binom{2m-1}{2k-1}\beta_k^s\right)\\
    &= \frac{1}{2^{2m-1}}\sum_{s=0}^{2m-2-l}(-1)^s\sum_{k=1}^m\left(\binom{2m-1}{2k-2}\alpha_k^s-\binom{2m-1}{2k-1}\beta_k^s\right)\\
    & =  \frac{1}{2^{2m-1}}\sum_{s=0}^{2m-2-l}(-1)^s c^s \sum_{k=1}^m\left(\binom{2m-1}{2k-2}(2k-2)^s-\binom{2m-1}{2k-1}(2k-1)^s\right)\\
    &=\frac{1}{2^{2m-1}}\sum_{s=0}^{2m-2-l}(-1)^s c^s \sum_{k=0}^{2m-1}\binom{2m-1}{k}(-1)^k k^s.
\end{align*}
Thus, to show $\mu_{0}(a)=\tilde{\mu}_0(a)$, it suffices to prove 
\begin{equation}
    \label{combination}
    \sum_{k=0}^{2m-1}\binom{2m-1}{k}(-1)^k k^s=0
\end{equation}
for all $s\leq 2m-2$. We will prove this by induction with respect to $s$. For $s=0$ \eqref{combination} holds trivially. Now suppose \eqref{combination} holds for $s$, we will prove that it also holds for $s+1$. Consider the generating function 
$$g(x) = (1+x)^{2m-1} = \sum_{k=1}^{2m-1}\binom{2m-1}{k}x^k.$$
Taking $s+1$-th order derivatives on both sides of the equation to obtain
$$C_{m,s}(1+x)^{2m-1-s} = \sum_{k=1}^{2m-1}\binom{2m-1}{k}\prod_{j=0}^{s}(k-j)x^{k-s+1}.$$
Now let $x=-1$, using the induction hypothesis, we have
$$0=(-1)^{1-s}\sum_{k=1}^{2m-1}(-1)^k \prod_{j=0}^s(k-j) = (-1)^{1-s}\sum_{k=1}^{2m-1}(-1)^k k^{s+1}.$$
This proves \eqref{combination}, thus $\mu_0=\tilde{\mu}_0$. We are done.
\end{proof}
\section{Proof of Theorem \ref{Angle_thm}}
\label{Proof_angle}
 We will first introduce some technical lemmas. 
    \begin{lemma}
    \label{lemma_test}
    Let $f_1,\dots ,f_m\in L^2(\real)$ be density functions such that $\|f_{k}\|_2\geq C_0$ for all $k=1,2,\cdots,m$ and $C_0>0$. Suppose $\tilde{f}\in L^2(\real)$ is a density function such that
    \begin{equation}
    \label{eq:angle-condition}
    |\langle \tilde{f},f_k\rangle|\leq \delta \|\tilde{f}\|_2\|f_k\|_2  \text{ for all $k\in[m]$ with $\delta<1$}.
    \end{equation}
     Then there exists a test function $\|w\|_2=1$, such that for all $k\in[m]$,
     $$\langle w,\tilde{f}\rangle=0,|\langle w,f_k\rangle|\geq  \frac{C_0\sqrt{1-\delta^2}}{4m^{3/2}}.$$
\end{lemma}
\begin{proof}
    Suppose $\calV\triangleq \Span\{\tilde{f},f_1,\dots ,f_m\}$ has dimension $r$. 
    Let $h_0 = \tilde{f}/\|\tilde{f}\|_2$, and let $h_1\dots ,h_{r-1}$ be an orthonormal basis for the orthogonal complement of $\Span\{\tilde f\}$ within $\calV$.
    Write $\tilde{f},f_1,\dots ,f_m$ as linear combinations of the orthonormal basis $h_0,h_1\dots ,h_{r-1}$:
    \begin{align*}
        \tilde{f}& =  \tilde{a}_0 h_0 + \sum_{i=1}^{r-1} \tilde{a}_i h_i,\\
        \tilde{f}_k& =  a_{k,0} h_0 + \sum_{i=1}^{r-1} a_{k,i} h_i,\quad k=1,\dots,m,
    \end{align*}
    where $\tilde{a}_0=\|\tilde{f}\|_2>0$ and $ \tilde{a}_i=0$ for $i=1,\ldots,r-1$. 
    It follows from the condition~\eqref{eq:angle-condition} that
    \begin{align*}
    &(\tilde{a}_0 a_{k,0})^2 \le \delta^2 \pth{\tilde{a}_0^2}
    \pth{a_{k,0}^2 + \sum_{i=1}^{r-1} a_{k,i}^2}\\
    \implies & \sum_{i=1}^{r-1} a_{k,i}^2 \ge (1-\delta^2)\sum_{i=0}^{r-1} a_{k,i}^2=(1-\delta^2)\|f_k\|_2^2.
    \end{align*}

    We then prove the lemma by the probabilistic method. 
    Let $t_1,\dots ,t_{r-1}\iiddistr \mathcal{N}(0,1)$ and $w'=\sum_{i=1}^{r-1}t_ih_i$.
    It suffices to show that the normalized function $w=w' / \|w'\|_2$ satisfies the desired property with strictly positive probability.
    By definition, $\|w\|_2=1$ and $\langle w,\tilde{f}\rangle=0$. 
    For a fixed $k\in[m]$, $\langle w',f_k\rangle = \sum_{i=1}^{r-1} a_{k,i} t_i \sim \mathcal{N}(0,\sum_{i=1}^{r-1} a_{k,i}^2)$. Let $\sigma_k^2 \triangleq \sum_{i=1}^{r-1} a_{k,i}^2$. 
    Then, 
    \[
    \prob\left[|\langle w',f_k\rangle|\leq \frac{\sqrt{2\pi}\sigma_k}{4m}\right]
    =2\prob\left[0\leq Z\leq \frac{\sqrt{2\pi}}{4m}\right] =2\int_{0}^{\frac{\sqrt{2\pi}}{4m}}\frac{1}{\sqrt{2\pi}}\exp\pth{-\frac{x^2}{2}}\diff x\leq \frac{1}{2m},
    \]
    where $Z$ is a standard Gaussian variable.
    Applying the union bound yields that
    \[
    \prob\left[|\langle w',f_k\rangle|\geq \frac{\sqrt{2\pi}\sigma_k}{4m},\forall k\in[m]\right]
    \geq 1-m\cdot\frac{1}{2m}=\frac{1}{2}.
    \]
    Moreover, since $\|w'\|_2^2\sim \chi_{r-1}^2$, by Markov inequality, we have $\prob[\|w'\|_2^2>4(r-1)]\le \frac{1}{4}$. Equivalently, $\prob[\|w'\|_2^2\le 4(r-1)]\ge \frac{3}{4}$.
    By the union bound, with probability at least $1/4$, 
    \[
    |\langle w,f_k\rangle|
    \ge \frac{1}{2\sqrt{r-1}}\cdot \frac{\sqrt{2\pi}\sigma_k}{4m}\geq \frac{C_0\sqrt{1-\delta^2}}{4m^{3/2}},\quad \forall k\in [m].
    \]
    This completes the proof.
\end{proof}

For the quantitative rates, we follow the concept of Kruskal rank and define the corresponding eigenvalues for a Gram matrix as follows.
\begin{definition}[Kruskal eigenvalue of a Gram matrix]
    Let $A\in\real^{m\times m}$ be a Gram matrix with. For $k\in [m]$, the $k$-th Kruskal eigenvalue of $A$ is defined as:
    $$\lambda_{k}^{\kru}(A):=\min\{ \lambda_{k}(A_{S\times S}): S\subseteq [m], |S|=k \},$$
    where $A_{S\times S}\in\real^{k\times k}$ is the principal submatrix of $A$ indexed by the set $S$.
\end{definition}
Evidently, 
if $h=(h_1,\dots,h_n)\in\mathcal{H}^n$ and $G=(\langle h_i,h_j\rangle)_{i,j=1}^n$ is the associated Gram matrix, then
$\lambda_{k}^{\kru}(G)>0$ implies $k_h\geq k$. We now present a lemma that establishes a lower bound for the Kruskal eigenvalue of the Hadamard product of two Gram matrices.
\begin{lemma}
    \label{condition_num_hadamard}
   Suppose $A,B\in\real^{m\times m}$ are Gram matrices.  Then for $k_1+k_2\leq m+1$, then
   $$\lambda^{\kru}_{k_1+k_2-1}(A\circ B)\geq \frac{\lambda_{k_1}^{\kru}(A)\lambda_{k_2}^{\kru}(B)}{k_1+k_2}.$$
\end{lemma}
\begin{proof}
    %
     Suppose $A=U^{\top}U,B=V^{\top}V$, where $U=[u_1,\dots,u_m]=A^{1/2}$, $V=[v_1,\dots,v_m]=B^{1/2}$. 
     Then $(A\circ B)_{ij}=(u_i^\top u_j)(v_i^\top v_j)=(u_i\otimes v_i)^\top (u_i\otimes v_i)$.
     Let $U\odot V=(u_{1}\otimes  v_1,\dots ,u_{m}\otimes v_m)$ denote the Khatri-Rao product.
     Then $A\circ B=(U\odot V)^\top (U\odot V)$.
     Consequently, 
     \begin{align*}
         \lambda_{k_1+k_2-1}^{\kru}(A\circ B) &= \min\{\lambda_{k_1+k_2-1}\left((A\circ B)_{S\times S}\right):S\subseteq [m], |S|=k_1+k_2-1\}\\
         &= \min\{\sigma^2_{k_1+k_2-1}((U\odot V)_{S}): S\subseteq [m], |S|=k_1+k_2-1\},
     \end{align*}
    where $(U\odot V)_S$ is the submatrix  containing the columns of $U\odot V$ indexed by $S$. Applying \cite[Lemma 20]{bhaskara_uniqueness_2013}, we have
    \begin{align*}
        &~\min\sth{\sigma^2_{k_1+k_2-1}((U\odot V)_{S}): S\subseteq [m], |S|=k_1+k_2-1}\\
        \geq&~\min\sth{\frac{\sigma_{k_1}^2(U_{S_1})\sigma_{k_2}^2(V_{S_2})}{k_1+k_2}:  |S_1|=k_1,|S_2|=k_2}\\
        \geq&~\frac{1}{k_1+k_2}\min\{\sigma_{k_1}^2(U_{S_1}):|S_1|=k_1\}\cdot \min\{\sigma_{k_2}^2(V_{S_2}):|S_2|=k_2\}\\
        =&~\frac{\lambda_{k_1}^{Kru}(A)\lambda_{k_2}^{Kru}(B)}{k_1+k_2}.
    \end{align*}
    The proof is completed.
\end{proof}

\begin{lemma}
\label{L_12_error}
Consider a Hilbert space $\calH=L^2(\Omega,\calF,\mu)$ with $\mu(\Omega)=1$. Let $f\in \calH$ satisfy $\|f\|_\infty \le C \|f\|_2$.
Suppose $g\in\calH$ and $\sin\theta (f,g) \leq \min\{\frac{\sqrt{3}}{2},\frac{1}{4C}\}$. 
Then
\[
\left\|\frac{f}{\|f\|_1}-\frac{g}{\|g\|_1}\right\|_2\leq 8C^2 \sin\theta (f,g).
\]
\end{lemma}
\begin{proof}
Without loss of generality, assume $\|f\|_2=\|g\|_2=1$.
Let $\theta= \theta (f,g)$. 
We decompose $g$ along $f$ and its orthogonal complement as 
\[
g= \cos \theta \cdot f+\sin \theta \cdot f^{\perp},
\]
where $\langle f,f^{\perp}\rangle=0$ and $\Norm{f^{\perp}}_2=1$.
Then, $ \Norm{g}_1 f-  \Norm{f}_1 g = (\Norm{g}_1- \Norm{f}_1 \cos\theta ) f - (\Norm{f}_1\sin\theta)  f^{\perp}.$
We obtain
\[
\left\|\frac{f}{\|f\|_1}-\frac{g}{\|g\|_1}\right\|_2
= \frac{\|\Norm{g}_1 f-  \Norm{f}_1 g\|_2}{\Norm{f}_1\Norm{g}_1}
= \frac{\sqrt{(\Norm{g}_1- \Norm{f}_1 \cos\theta )^2+(\Norm{f}_1\sin\theta)^2}}{\Norm{f}_1\Norm{g}_1}.
\]
By triangle inequality, $|\Norm{g}_1- \Norm{f}_1 \cos\theta|\le \Norm{g-f \cos\theta}_1=\Norm{f^{\perp}}_1 \sin\theta$.
By Cauchy-Schwarz inequality, $\Norm{f}_1\le \Norm{f}_2\le 1$ and $\Norm{f^{\perp}}_1\le \Norm{f^{\perp}}_2\le 1$.
It follows that
\begin{equation}
\label{eq:diff-L1-normalized}
\left\|\frac{f}{\|f\|_1}-\frac{g}{\|g\|_1}\right\|_2
\le \frac{\sqrt{2}\sin\theta}{\Norm{f}_1\Norm{g}_1}.
\end{equation}

It remains to lower bound $\Norm{f}_1$ and $\Norm{g}_1$.
Since $\|f\|_{\infty}\leq C$, we have
\[
1=\int f^2\diff \mu
\le C \int |f|\diff \mu
= C \Norm{f}_1.
\]
Furthermore, by the triangle inequality,
\[
\Norm{g}_1
\ge \cos \theta \Norm{f}_1 - \sin \theta \Norm{f^{\perp}}_1
\ge \frac{\cos \theta}{C} - \sin\theta.
\]
Since $\sin\theta\le \min\{\frac{\sqrt{3}}{2},\frac{1}{4C}\}$, we have $\frac{\cos \theta}{C} - \sin\theta\ge \frac{1}{4C}$.
The conclusion follows from \eqref{eq:diff-L1-normalized}.
\end{proof}

    \begin{lemma}[\cite{Gohberg1990} Corollary 1.6]
        \label{Weyl}
        Suppose $\mathcal{H}_1,\mathcal{H}_2$ are two Hilbert spaces, and $A,B:\mathcal{H}_1\mapsto \mathcal{H}_2$ are two finite rank operators with rank $\leq m$. Denote the singular values of $A,B$ by $\sigma_1(A)\geq\cdots\geq \sigma_{m}(A)\geq 0$ and $\sigma_1(B)\geq\cdots\geq \sigma_{m}(B)\geq 0$, respectively. Then we have
        \begin{equation*}
            \max_{k\in[m]}|\sigma_k(A)-\sigma_{k}(B)| \leq \|A-B\|_{\textnormal{op}}.
        \end{equation*}
    \end{lemma}
Now we are ready to prove Theorem \ref{Angle_thm}. 
\begin{proof}[Proof of Theorem \ref{Angle_thm}]
For $I\subseteq[d]$, let $f_I$ and $\tilde{f}_I$ denote the marginal densities of $f$ and $\tilde{f}$ with respect to the variables indexed by $I$, respectively. Let $x_I=(x_i)_{i\in I}\in [0,1]^{|I|}$ and $x_{-I}=(x_i)_{i\in I^c}\in [0,1]^{d-|I|}$. 
From Cauchy-Schwarz inequality, we have
\begin{align*}
    \|f_I-\tilde{f}_I\|_2 &= \int_{[0,1]^{d-|I|}}  \left(\int_{[0,1]^{|I|}}1\cdot \left( f(x_I,x_{-I})-\tilde{f}(x_I,x_{-I})\right)dx_I\right)^2 dx_{-I}\\
    &\leq \int_{[0,1]^{d-|I|}} \left(f(x_I,x_{-I})-\tilde{f}(x_I,x_{-I})\right)^2dx_Idx_{-I}\\
    &=\|f-\tilde{f}\|_2\leq \epsilon.
\end{align*}
Thus, we only need to prove the result for $d=2m-1$. 

We begin with some preliminary preparations. From $f_{kj},\tilde{f}_{kj}\leq C$, we know $f_{kj},\tilde{f}_{kj}\in L^2([0,1])$ for every $k\in[m]$ and $j\in [2m-1]$.
Thus, applying a unitary transformation $U$, we map $f,\tilde{f}$ to $T,\tilde{T}\in L^2([0,1])^{\otimes(2m-1)}$, respectively, with the following explicit form: 
$$T=\sum_{k=1}^m\pi_k\otimes_{j=1}^{2m-1}f_{kj},\ 
\quad
\tilde{T}=\sum_{k=1}^m\tilde{\pi}_k\otimes_{j=1}^{2m-1}\tilde{f}_{kj}.$$ 
We consider the following transform:
For $w\in L^2([0,1])$, we write the mode-$1$ multiplication of $T$ as 
\begin{equation}
    \label{eq:mode1}
    T\times _1w = \sum_{k=1}^m \pi_k\langle w,f_{k1}\rangle\otimes_{j=2}^{2m-1}f_{kj} \in L^2([0,1])^{\otimes 2m-2}.
\end{equation}
Then, applying a unitary transformation $U'$, we unfold $T\times_1 w$ to the following linear operator: 
\begin{equation}
    \label{eq:key_mat}
    T_{w}=AD_{\pi,w}B^* \in L^2([0,1])^{\otimes (m-1)}\otimes L^2([0,1])^{\otimes (m-1)},
\end{equation}
where $A=(\otimes_{j=2}^m f_{1j},\dots ,\otimes_{j=2}^m f_{mj}), B=(\otimes_{j=m+1}^{2m-1} f_{1j},\dots ,\otimes_{j=m+1}^{2m-1} f_{mj}),D_{\pi,w}=\diag\{\pi_1$\\
$\langle w,f_{11}\rangle,\dots ,\pi_m\langle w,f_{m1}\rangle\}$, and $B^*$ is the adjoint operator of $B$. 
Similarly, we map $\tilde{T}$ to $\tilde{T}_w=\tilde{A}D_{\tilde{\pi},w}\tilde{B}^*$.
Note that $U, U'$ are both unitary and therefore preserves the inner product, we deduce that $\|T-\tilde{T}\|_{\op}=\|f-\tilde{f}\|_2\leq \epsilon$, $\|T_w-\tilde{T}_w\|_{\op}=\|T\times_1 w-\tilde{T}\times_1w\|_{\op}$.
Additionally, we have the following relation:
 \begin{align}
 \label{eq:opnorm}
        \sup_{w\in L^2([0,1]),\|w\|_2=1}\|T\times_1w-\tilde{T}\times_1w\|_{\textnormal{op}} &= \sup_{\substack{w\in L^2([0,1]),\|w\|_2=1 \nonumber\\ \tilde{w}\in L^2([0,1]^{2m-2}),\|\tilde{w}\|_2=1}} 
\langle T\times_1w-\tilde{T}\times_1w, \tilde{w}\rangle\nonumber\\ 
        &= \sup_{\substack{w\in L^2([0,1]),\|w\|_2=1 \nonumber\\ \tilde{w}\in L^2([0,1]^{2m-2}),\|\tilde{w}\|_2=1}} 
\langle T-\tilde{T}, w\otimes\tilde{w}\rangle\nonumber\\
        &\leq \sup_{w'\in L^2([0,1])^{\otimes (2m-1)},\|w'\|_2=1} \langle T-\tilde{T},w'\rangle\nonumber\\
        &=\|T-\tilde{T}\|_{\textnormal{op}}\leq \epsilon.
    \end{align}
Thus, $\sup_{\|w\|_2=1}\|T_w-\tilde{T}_w\|_{\op}\leq \epsilon$.
Note that $T_w,\tilde{T}_w$ are both finite rank linear operators with rank at most $m$. By Lemma \ref{Weyl}, we have
\begin{equation}
    \label{eq:singular_diff}
    \sup_{w\in L^2([0,1]),\|w\|_2=1}\max_{k\in[m]}|\sigma_{k}(T_w)-\sigma_k(\tilde{T}_w)|\leq \epsilon.
\end{equation}
Now we show that in \eqref{eq:key_mat}, $A, B$ are well conditioned as finite rank linear operators, which allows us to focus on the diagonal matrix $D_{w,\pi}$ afterwards. Iteratively applying Lemma \ref{condition_num_hadamard} with $k_1=2$, we have a lower bound of the $m$-th singular value of $A$:
    \begin{equation}
        \label{condition_num}
        \sigma_m(A)=\sqrt{\lambda_{m}^{Kru}(A^*A)} =\sqrt{\lambda_{m}^{\kru}(A_2\circ A_3\circ\dots\circ A_m)}\geq \sqrt{\frac{\prod_{j=2}^m\lambda_2^{\kru}(A_j)}{(m-1)!}}\geq \sqrt{\frac{(1-\mu)^{m-1}}{(m-1)!}},
    \end{equation}
    where $A_j$ is the Gram matrix of $f_j=(f_{1j},...,f_{mj})$. The last inequality is because $\|f_{kj}\|_2\geq 1$ and the incoherence condition. Similarly, $\sigma_m(B)\geq \sqrt{\frac{(1-\mu)^{m-1}}{(m-1)!}}$. 
    
We prove Theorem \ref{Angle_thm} by contradiction, showing that it conflicts with equation \eqref{eq:singular_diff} for some $\|w\|_2=1$ and $k\in[m]$. The proof is divided into the following four steps.

 \textbf{Step 1: Find a component density close to the true one}: Define $\epsilon'\triangleq \frac{L_m}{(1-\mu)^{m-1}\zeta}\epsilon$.  We show that for any $(k,j)\in [m]\times[2m-1]$, there exists $k'\in[m]$ such that $\|f_{k'j}-\tilde{f}_{kj}\|_2\le 8C^2 \epsilon'$ for every $j\in [2m-1]$; Without loss of generality, we show this for $j=1$. 
 From Cauchy-Schwarz inequality, we have $\|f_{k'1}\|_2\geq\|f_{k'1}\|_1=1$ and thus $\|f_{k'1}\|_{\infty}\leq C \|f_{k'1}\|_2$. 
 From the assumption on $\epsilon$, we can verify $\epsilon'\leq \frac{1}{4C}\wedge \frac{\sqrt{3}}{2}$.
 Thus, by Lemma \ref{L_12_error}, it suffices to show $\sin\theta(f_{k'1},\tilde{f}_{k1})\le \epsilon'$.
 
Suppose on the contrary there exists some $k\in[m]$ such that for all $k'\in[m]$, $\sin\theta(f_{k'1},\tilde{f}_{k1})> \epsilon'.$
Consequently, $|\langle f_{k'1},\tilde{f}_{k1}\rangle|\leq \sqrt{1-\epsilon'^2}\|f_{k'1}\|_2\|\tilde{f}_{k1}\|_2$ for all $k'\in[m]$.
By Lemma \ref{lemma_test}, there exists $w_0\in L^2([0,1])$ with $\|w_0\|_2=1$ such that 
    \begin{equation}
        \forall k'\in[m],\ |\pi_k\langle w_0,f_{k'1}\rangle|\geq \frac{\zeta\epsilon'}{4m^{3/2}}=\frac{(m-1)!}{(1-\mu)^{m-1}}\epsilon,\ \langle w_0,\tilde{f}_{k1}\rangle=0.
    \end{equation}
Thus, the diagonal matrix $D_{\tilde{\pi},w_0}$ has a zero diagonal entry, which implies that $\sigma_{m}(\tilde{T}_{w_0})=0$.
On the other hand, 
\begin{equation}
    \label{eq:singular_of_D}
    |\sigma_m(D_{\pi,w_0})|\geq \min_{k}|\pi_k\langle w_0,f_{k1}\rangle|\geq\frac{(m-1)!}{(1-\mu)^{m-1}}\epsilon.
\end{equation}
Thus, we obtain
\begin{equation*}
    |\sigma_m(T_{w_0})-\sigma_m(\tilde{T}_{w_0})|=|\sigma_{m}(T_{w_0})|\geq \sigma_m(A)\sigma_m(B)|\sigma_{m}(D_{\pi,w})|>\epsilon,
\end{equation*}
a contradiction to \eqref{eq:singular_diff}.
    \textbf{Step 2: Verify the mapping is one-to-one}. We will show that the mapping $\sigma_j: k\mapsto k'$ in Step 1 is one-to-one for every $j\in[2m-1]$, thus a permutation.  
    Suppose this is not true, then there exists $j\in[2m-1]$ and $k_1,k_2,k'\in[m], k_1\ne k_2$, such that $\|\tilde{f}_{k_1j}-f_{k'j}\|_2,\|\tilde{f}_{k_2 j}-f_{k'j}\|_2\leq 8C^2\epsilon'$. Without loss of generality, take $k'=j=1$. 
    For $f_{21},...,f_{m1}$ and $f_{11}$ $\mu$-incoherent with them, applying Lemma \ref{lemma_test}, there exists a $w_1\in L^2(\real)$ with $\|w_1\|_2=1$, such that
    \begin{equation}
        \label{w_1}
        \langle w_1,f_{11}\rangle=0, |\langle w_1,f_{k1}\rangle|\geq \frac{\sqrt{1-\mu^2}}{4m^{3/2}},\ k=2,3,...,m.
    \end{equation}
    Since $\|f_{11}-\tilde{f}_{k_11}\|_2\leq\epsilon'$, we know 
    $$|\langle w_1, \tilde{f}_{k_11}\rangle|=|\langle w_1, \tilde{f}_{k_11}-f_{11}\rangle|\leq \|\tilde{f}_{k_11}-f_{11}\|_2\leq 8C^2\epsilon'.$$
    Similarly, $|\langle w_1,\tilde{f}_{k_21}\rangle|\leq 8C^2\epsilon'$. 
    Consequently, $\sigma_{m-1}(D_{\pi,w_1})\geq \frac{\zeta\sqrt{1-\mu^2}}{4m^{3/2}}$, whereas $|\sigma_{m-1}(D_{\tilde{\pi},w_1})|\leq \epsilon'$. 
    Similar to Step 1, we deduce that
    \begin{align*}
        |\sigma_{m-1}(T_{w_1})-\sigma_{m-1}(\tilde{T}_{w_1})|&\geq|\sigma_{m-1}(T_{w_1})| - |\sigma_{m-1}(\tilde{T}_{w_1})|\\
        &\geq \sigma_m(A)\sigma_m(B)|\sigma_{m-1}(D_{\pi,w})|-\sigma_1(\tilde{A})\sigma_1(\tilde{B})|\sigma_{m-1}(D_{\tilde{\pi},w_1})|\\
        &\geq \frac{(1-\mu)^{m-1}}{(m-1)!}\frac{\zeta\sqrt{1-\mu^2}}{4m^{3/2}} - C^{2m-2}\epsilon'\\
        &= \frac{(1-\mu)^{m-1}\zeta\sqrt{1-\mu^2}}{L_m}-\frac{8C^{2m}L_m}{(1-\mu)^{m-1}}\epsilon>\epsilon,
    \end{align*}
 a contradiction to \eqref{eq:singular_diff}. The last inequality is from the assumption on $\epsilon$. This proves that $\sigma_j$ is an injection from $[m]$ to $[m]$, thus a permutation.
\par \textbf{Step 3: Show that the permutations are identical.} We will prove that $\sigma_1=\cdots =\sigma_{2m-1}$.
Suppose on the contrary there exists $j_1,j_2\in [2m-1]$ such that $\sigma_{j_1}\ne \sigma_{j_2}$; without loss of generality, we take $j_1=1,j_2=2$. From $\sigma_{1}\ne\sigma_2$, there exists $k_1,k_2\in [m], k_1\ne k_2$ such that $\sigma_1(k_1)=\sigma_2(k_2)$; without loss of generality, we take $\sigma_{1}(1)=\sigma_2(2)=1$.
From the triangle inequality, we have
\begin{align*}
    &\quad\left\|\sum_{k=1}^m(f_{\sigma_1(k)1}-\tilde{f}_{k1})\otimes\tilde{f}_{k2}\otimes(\tilde{\pi}_k\otimes_{j=3}^{2m-1}\tilde{f}_{kj})\right\|_{\op}\\
    &\leq \sum_{k=1}^m\left\|(f_{\sigma_1(k)1}-\tilde{f}_{k1})\otimes\tilde{f}_{k2}\otimes(\tilde{\pi}_k\otimes_{j=3}^{2m-1}\tilde{f}_{kj})\right\|_{\op} \\
    &\leq m\cdot8C^2\epsilon' \cdot C^{2m-2}=8mC^{2m}\epsilon'.
\end{align*}
Similarly,
$$\left\|\sum_{k=1}^m f_{\sigma_{1}(k)1}\otimes(f_{\sigma_2(k)2}-\tilde{f}_{k2})\otimes(\tilde{\pi}_k\otimes_{j=3}^{2m-1}\tilde{f}_{kj})\right\|_{\op}\leq 8mC^{2m}\epsilon'.$$
Let $T'=\sum_{k=1}^m\tilde{\pi}_k f_{\sigma_{1}(k)1}\otimes f_{\sigma_2(k)2}\otimes(\otimes_{j=3}^{2m-1}\tilde{f}_{kj})$. From the triangle inequality and $\|T-\tilde{T}\|_{\op}\leq \epsilon$, we deduce that
\begin{align*}
   \|T-T'\|_{\op}= &\ \left\|\sum_{k=1}^{m}\left( f_{k1}\otimes f_{k2}\otimes (\pi_k\otimes_{j=3}^{2m-1}f_{kj})- f_{\sigma_{1}(k)1}\otimes f_{\sigma_2(k)2}\otimes(\tilde{\pi}_k\otimes_{j=3}^{2m-1}\tilde{f}_{kj})\right)\right\|_{\op}\\
    \leq&\ \left\|\sum_{k=1}^m(f_{\sigma_1(k)1}-\tilde{f}_{k1})\otimes\tilde{f}_{k2}\otimes(\tilde{\pi}_k\otimes_{j=3}^{2m-1}\tilde{f}_{kj})\right\|_{\op}\\
    +&\ \left\|\sum_{k=1}^m f_{\sigma_{1}(k)1}\otimes(f_{\sigma_2(k)2}-\tilde{f}_{k2})\otimes(\tilde{\pi}_k\otimes_{j=3}^{2m-1}\tilde{f}_{kj})\right\|_{\op}+\|T-\tilde{T}\|_{\op}\\
    \leq &\ \epsilon+16m C^{2m}\epsilon'.
\end{align*}
Since $\{f_{k1}\}_{k=1}^m,\{f_{k2}\}_{k=1}^m$ are $\mu$-incoherent, by applying Lemma \ref{lemma_test} again, there exists $u,v\in L^2([0,1])$ with $\|u\|_2=\|v\|_2=1$ such that 
$$\langle u,f_{11}\rangle=\langle v,f_{12}\rangle=0,|\langle u,f_{k1}\rangle|,|\langle v,f_{k2}\rangle|\geq \frac{\sqrt{1-\mu^2}}{4m^{3/2}},\ k=2,3,\cdots m.$$
Let  $\times_j$ denote the mode-$j$ multiplication of a tensor. For $w\in L^{2}([0,1])$, define $T_{u,v,w}\triangleq T\times_1u\times_2v \times_3 w, T_{u,v,w}'\triangleq T'\times_1u\times_2v\times_3 w$, respectively. 
Then $T_{u,v,w},T'_{u,v,w}\in L^2([0,1])^{\otimes (2m-4)}$.
From $\sigma_1(1)=\sigma_2(2)=1$ and the choice of $u,v$, we obtain
$$T_{u,v,w}=\sum_{k=2}^m\langle f_{k1},u\rangle\langle f_{k2},v\rangle \langle f_{k3},w\rangle \pi_k\otimes_{j=4}^{2m-1}f_{kj},$$
and
$$T'_{u,v,w}=\sum_{k=3}^m\langle f_{\sigma_{1}(k)1},u\rangle \langle f_{\sigma_2(k)2},v\rangle \langle\tilde{f}_{k3},w\rangle\tilde{\pi}_k\otimes_{j=4}^{2m-1}\tilde{f}_{kj}.$$
By applying a unitary transform, we unfold $T_{u,v,w}$ to 
$$S_{u,v,w} = A_1 D_{u,v,w,\pi}B_1^*\in L^2([0,1])^{\otimes (m-2)}\otimes L^2([0,1])^{\otimes (m-2)},$$
where $A_1 = (\otimes_{j=4}^{m+3}f_{2j},\cdots,\otimes_{j=4}^{m+3}f_{2j}), B_1=(\otimes_{j=m+4}^{2m-1}f_{2j},\cdots,\otimes_{j=m+4}^{2m-1}f_{2j})$, and 
\newline $D_{u,v,w,\pi}= \diag \left(\pi_2\langle f_{21},u\rangle\langle f_{22},v\rangle\langle f_{23},w\rangle,\cdots \pi_m\langle f_{m1},u\rangle\langle f_{m2},v\rangle\langle f_{m3},w\rangle \right)$.
Similarly, denote the image of $T'_{u,v,w}$ by $S'_{u,v,w}=\tilde{A}_1 D_{u,v,w,\tilde{\pi}}\tilde{B}_1^*$.
Similar to \eqref{eq:singular_diff}, we have $$\sup_{\|w\|_2=1}\max_{k\in [m]}\|\sigma_{k}(S_{u,v,w})-\sigma_k(S'_{u,v,w})\|_{\op}\leq \|T-T'\|_{\op}\leq \epsilon+16mC^{2m}\epsilon'<17mC^{2m}\epsilon'.$$
From Lemma \ref{condition_num_hadamard} again, we have $\sigma_{m-1}(A_1),\sigma_{m-1}(B_1)\geq \sqrt{\frac{(1-\mu)^{m-2}}{(m-2)!}}$.
Since $T'_{u,v,w}$ has rank at most $m-2$, we obtain $\sigma_{m-1}(S'_{u,v,w})=0$ for any $w$.
Thus, choosing $w=\frac{f_{23}}{\|f_{23}\|_2}$, we obtain
\begin{align*}
    |\sigma_{m-1}(S_{u,v,w})-\sigma_{m-1}(S'_{u,v,w})|&=|\sigma_{m-1}(S_{u,v,w})|\geq \sigma_{m-1}(A)\sigma_{m-1}(B)|\sigma_{m-1}(D_{u,v,w})|\\
    &\geq \frac{(1-\mu)^{m-2}}{(m-2)!}\cdot \frac{\zeta(1-\mu^2)(1-\mu)}{16m^3}\\
    &>17mC^{2m}\epsilon',
\end{align*}
which leads to a contradiction. The last inequality follows from the assumption on $\epsilon$.
This proves $\sigma_j$'s are identical.
    \par\textbf{Step 4: Bounding the error of mixing proportion}. For the remainder of this proof, we assume $\sigma$ is the identity without loss of generality. In this step, the norm $\|\cdot\|$ refers to the operator norm if not specified. We consider the marginal density on the first $m-1$ variables:
    \begin{equation}
        \label{marginal_joint}
        f_{1:(m-1)} = \sum_{k=1}^m \pi_k \prod_{j=1}^{m-1} f_{kj} = F_1 \pi,
    \end{equation}
    where $F_1=(\prod_{j=1}^{m-1}f_{1j},\dots ,\prod_{j=1}^{m-1}f_{mj})$, a rank-$m$ linear operator from $\real^m$ to $L^2(\real^{m-1})$. Similarly, we define $\tilde{f}_{1:(m-1)},\tilde{F}_1 $ and $\tilde{\pi}$ from $\tilde{f}$.
    Let $\tilde{f}_{1:(m-1)} - f_{1:(m-1)}= h, \tilde{F}_1 -F_1 = E_2$, and $\pi -\tilde{\pi} = e_3$. We have
    \begin{align*}
        \tilde{f}_{1:(m-1)} = \tilde{F}_1 \tilde{\pi}&\implies (f_{1:(m-1)}+h) = (F_1+E_2)(\pi+e_3)\\
        &\implies\tilde{F}_1 e_3 = h - E_2\pi.
    \end{align*}
    Since $F_1,\tilde{F}_1$ are both rank-$m$, by Lemma \ref{Weyl},  $\sigma_{m}(\tilde{F}_1) \geq \sigma_{m}(F_1) - \|E_2\|$. 

     Now we give an upper bound of $\|E_2\|$. We first bound $\sin\theta(\prod_{j=1}^{m-1}f_{kj},\prod_{j=1}^{m-1}\tilde{f}_{kj})$:
    \begin{align*}
        \label{cal_sin}
        \sin\theta\left(\prod_{j=1}^{m-1}f_{kj},\prod_{j=1}^{m-1}\tilde{f}_{kj}\right)&= \sqrt{1-\cos^2\theta\left(\prod_{j=1}^{m-1}f_{kj},\prod_{j=1}^{m-1}\tilde{f}_{kj}\right)}\\
        & = \sqrt{1-\prod_{j=1}^{m-1}\cos^2\theta\left(f_{kj},\tilde{f}_{kj}\right)}\\
        &\leq \sqrt{1-(1-\epsilon'^2)^{m-1}} \\
        &\leq \epsilon'\sqrt{m-1}. 
    \end{align*}
We have $\left\|\prod_{j=1}^{m-1}f_{kj}\right\|_\infty \le C^{m-1}\le C^{m-1}\left\|\prod_{j=1}^{m-1}f_{kj}\right\|_2$,
by Lemma \ref{L_12_error}, we have 
    \begin{equation*}
        \left\|\prod_{j=1}^{m-1}f_{kj}-\prod_{j=1}^{m-1}\tilde{f}_{kj}\right\|_2\leq 8C^{2m-2}\sqrt{m-1}\epsilon'.
    \end{equation*}
    Thus, 
    \begin{align}
    \label{eq:Upper_E2}
        \|E_2\| &= \|F_1-\tilde{F}_1\| = \sup_{\|x\|_2=1}\|(F_1-\tilde{F}_1)x\|_2\nonumber\\
        &= \sup_{\|x\|_2=1}\left\|\sum_{k=1}^m x_k(\prod_{j=1}^{m-1}f_{kj}-\prod_{j=1}^{m-1}\tilde{f}_{kj})\right\|_2\nonumber\\
        &\leq\sum_{k=1}^m \left\|\prod_{j=1}^{m-1}f_{kj}-\prod_{j=1}^{m-1}\tilde{f}_{kj}\right\|_2\nonumber\\
        &\leq 8C^{2m-2}m\sqrt{m-1}\epsilon'.
    \end{align}
    From the assmuption on $\epsilon$, we know $\|E_2\|_2\leq \frac{1}{2}\sqrt{\frac{(1-\mu)^{m-1}}{(m-1)!}}\leq \frac{1}{2}\sigma_{m}(F_1)$, thus $\sigma_{m}(\tilde{F}_1)\geq \frac{1}{2}\sigma_{m}(F_1)>0$.
    From the triangle inequality,
    \begin{equation*}
        \|h-E_2\pi\|_2\leq \|h\|_2+\|E_2\|\|\pi\|_2 \leq \epsilon+\|E_2\|.
    \end{equation*}
    Thus, plugging in \eqref{eq:Upper_E2}, we obtain the upper bound of $\|\pi-\tilde{\pi}\|_2$:
    $$\|\pi-\tilde{\pi}\|_2=\|e_3\|_2 =\| \tilde{F}_1^{-1}(h-E_2\pi)\|_2 \leq \frac{2}{\sigma_{m}(F_1)}\cdot(\epsilon+\|E_2\|) \leq \frac{16C^{2m-2}L_m^2}{(1-\mu)^{\frac{3(m-1)}{2}}\zeta}\epsilon$$
    as desired.
\end{proof}

\section{Proof of Theorem \ref{Holder_rate}}
\label{pf:stat_limit}
\subsection{Definitions and some preparations}
\label{entropy_minimax}
For the H\"{o}lder class, we give a formal definition for the H\"{o}lder smooth function in the main text:
\begin{definition}
\label{Holder}
    For a parameter $q=l+\beta>0$, where $l\in \mathbb{Z},\beta\in (0,1]$, we say a function $f$ is $q$-H\"{o}lder smooth with parameter $L>0$, if $f$ is $l$-times continuously differentiable, and the $l$-th derivative satisfies
    $$\left|\frac{\diff^l f}{\diff x^l}(x)-\frac{\diff^l f}{\diff x^l}(y)\right|\leq L|x-y|^{\beta}.$$
\end{definition}
 Now we review some classical results about metric entropy that we need for the proof. We begin from the concept of metric entropy.
\begin{definition}[covering and packing entropy]
    \label{cov_pac}
    Let $\mathcal{F}$ be a class of densities and $\rho$ be a metric. 
    \begin{enumerate}
        \item  An $\epsilon$-packing of $\calF$ with respect to $\rho$ is a subset $\calM = \{f_1,\dots ,f_M\}\subset \mathcal{F}$ such that $\rho(f_i,f_j)\geq \epsilon$ for all $i\neq j$. The $\epsilon$-packing number of $\mathcal{F}$ is defined to be the maximum number $M=M(\mathcal{F},\rho,\epsilon)$ such that there exists a $\epsilon$-packing with cardinality $M$.
        \item  An $\epsilon$-net of $\mathcal{F}$  with respect to $\rho$ is a set $\mathcal{N}=\{f_1,\dots ,f_N\}$ such that, for all $f\in\mathcal{F}$, there exists $i\in [N]$ such that $\rho(f_i,f)<\epsilon$. The $\epsilon$-covering number is defined to be the minimum $N=N(\mathcal{F},\rho,\epsilon)$ such that there exists a $\epsilon$-net with cardinality $N$.
    \end{enumerate}
    The $\epsilon$-covering entropy and $\epsilon$-packing entropy are defined as the logarithm of the $\epsilon$-covering number and $\epsilon$-packing number, respectively.
\end{definition}
For a class $\mathcal{F}$ and a metric $\rho$, there is a well-known relationship between covering and packing number \cite[see e.g.][Theorem 27.2]{Polyanskiy_Wu_2025}:
\begin{equation}
    \label{cov_pack}
    M(\mathcal{F},\rho,2\epsilon)\leq N(\mathcal{F},\rho,\epsilon)\leq M(\mathcal{F},\rho,\epsilon).
\end{equation}
There is a close relationship between the entropy of a class and the minimax risk. For the minimax upper bound, we have the following classical results from \cite{yatracos1985rates,birge1983approximation}:
\begin{proposition}
\label{results_mini_upper_classic}
     For $\rho\in \{\textnormal{TV},H\}$ and a class of density $\mathcal{F}$, given a random sample $X_1,\dots ,X_n\sim f\in\mathcal{F}$, we have entropic minimax upper bounds:
\begin{equation*}
 \inf_{\hat{f}}\sup_{f\in\mathcal{F}}\Expect[\rho^2(\hat{f},f)]\lesssim \inf_{\epsilon>0}\left\{\epsilon^2 + \frac{1}{n} \log N(\mathcal{F},\rho,\epsilon)\right\},
\end{equation*}
\end{proposition}
We can also derive the minimax lower bound from the bounds of metric entropy. The fundamental work of this characterization is from \cite{yang_information-theoretic_1999}.
\begin{proposition}[Theorem 1 in \cite{yang_information-theoretic_1999}]
\label{Yang_minimax}
Let $KL(f||g):=\int f(x)\log\frac{f(x)}{g(x)}\diff x$ be KL-divergence between $f$ and $g$. The KL $\epsilon$-covering number for a class of densities $\mathcal{F}$ is defined by 
    \begin{equation*}
        N(\mathcal{F},\sqrt{\KL},\epsilon):=\min\{N: \exists q_1,\dots q_N\ s.t.\forall f\in\mathcal{F},\exists i\in[N],\KL(f||q_i)\leq\epsilon^2\}.
    \end{equation*} Define the covering radius $\epsilon_n$ of $\mathcal{F}$ to be the solution of the following equation:
\begin{equation}
    \label{KL_radius}
    \epsilon_n^2 = N(\mathcal{F},\sqrt{\KL},\epsilon_n)/n.
\end{equation}
Suppose we are given a random sample $X_1,\dots, X_n\sim f\in\mathcal{F}$. Then, for any metric $\rho$ with triangle inequality, the minimax risk has a lower bound :
\begin{equation}
    \label{mini_lower}
    \inf_{\hat{f}}\sup_{f\in\mathcal{F}}\Expect[\rho^2(\hat{f},f)] \geq \frac{1}{8}\epsilon_{n,\rho}^2,
\end{equation}
where $\epsilon_{n,\rho}$ is defined by the equation
\begin{equation}
    \label{Packing_lower}
    M(\mathcal{F},\rho,\epsilon_{n,\rho}) = 4n\epsilon_n^2 +2\log 2.
\end{equation}
\end{proposition}
For calculating the cardinality of a packing set, we use the following result.
\begin{proposition}[Gilbert–Varshamov bound]
\label{GV_bound}
    Let $A_{M,n} = \{1,2,\dots,M\}^{n}$. For $a=(a_1,\dots,a_n),b=(b_1,\dots,b_n)\in A_{M,n}$, define the Hamming distance of $a,b$ to be
\begin{equation*}
    \textnormal{Ham}(a,b) = \|a-b\|_0:= \#\{i\in[n]: a_i\ne b_i\}.
\end{equation*}
Let $P_{M,n}(d)$ be a $d$-packing of $A_{M,n}$ with respect to Hamming distance. Then for $d\leq n$,
\begin{equation*}
    |P_{M,n}(d)| \geq \frac{M^n}{\sum_{j=0}^{d-1}\binom{n}{j}(M-1)^j}.
\end{equation*}
\end{proposition}
\subsection{Entropic bounds}
We will first prove the following entropic bounds. 
\begin{lemma}
\label{entro_lower_bound}
Let $\calF_{L,q}$ denote the class of all $q$-H\"{o}lder smooth densities on $[0,1]$ with smoothness parameter $q$ and constant $L>0$. Let $\calG_{\calF_{L,q}}$ be defined as in \eqref{MM_class}. Then we have
    $$ d\left(\frac{1}{\epsilon}\right)^{1/q}\lesssim_{L,q}\log N(\mathcal{G}_{\calF_{L,q}}^{(m,d)},\TV,\epsilon)\lesssim_{L,q}md^{1+\frac{1}{q}}\left(\frac{1}{\epsilon}\right)^{1/q}\quad \forall \epsilon>0.$$
$$ d^{1+\frac{1}{q}}\left(\frac{1}{\epsilon}\right)^{2/q}\lesssim_{L,q}\log N(\mathcal{G}_{\calF_{L,q}}^{(m,d)},H,\epsilon)\lesssim_{L,q} md^{1+\frac{1}{q}}\left(\frac{1}{\epsilon}\right)^{2/q}.\quad\forall 0<\epsilon<1.$$ 
\end{lemma}
\begin{proof}[Proof of Lemma \ref{entro_lower_bound}]
\textbf{Upper bound:} We first prove the entropic upper bound under $\TV$. Pick a $\epsilon/2d$-covering of $\mathcal{F}_{L,q}$ under $\TV$, denoted by $S=\{h_1,\dots ,h_{|S|}\}$. 
     Also, pick a $\epsilon/2$-covering of the simplex $\Delta^{m-1}$, denoted by $D_{\epsilon/2}$. We consider the following set:
    $$\mathcal{N}=\left\{\tilde{f}=\sum_{k=1}^{m}\tilde{\pi}_k\prod_{j=1}^d\tilde{f}_{kj}(x_j):\tilde{f}_{kj}\in S, \tilde{\pi}=(\tilde{\pi}_1,\dots ,\tilde{\pi}_m)\in D_{\epsilon/2}\right\}.$$
     We now prove that $\mathcal{N}$ is indeed an $\epsilon$-covering of $f\in \mathcal{G}_{L,q}^{(m,d)}$. For any $\mathcal{G}_{L,q}^{(m,d)}$, there exists an element in $\tilde{f}\in\mathcal{N}$ such that
    $$\tilde{f} = \sum_{k=1}^m\tilde{\pi}_k\prod_{j=1}^d\tilde{f}_{kj},\left\|f_{kj}-\tilde{f}_{kj}\right\|_1\leq \epsilon/d, \ \forall k,j; \left\|\pi-\tilde{\pi}\right\|_1\leq \epsilon/2.$$
    By triangle inequality, we have (all integrals are under Lebesgue measure)
    \begin{align}
        \left\|f-\tilde{f}\right\|_1&\leq \int\left|\sum_{k=1}^m\pi_k\prod_{j=1}^d f_{kj}-\sum_{k=1}^m\tilde{\pi}_k\prod_{j=1}^d\tilde{f}_{kj}\right|\nonumber\\
        &\leq \sum_{k=1}^m\int \left|\pi_k\prod_{j=1}^d f_{kj}-\tilde{\pi}_k\prod_{j=1}^d\tilde{f}_{kj}\right|\nonumber\\
        &\leq \sum_{k=1}^m\left(\int  \left|\pi_k-\tilde{\pi}_k\right|\prod_{j=1}^d f_{kj}
        +\tilde{\pi}_k\int \left|\prod_{j=1}^d f_{kj}-\prod_{j=1}^d\tilde{f}_{kj}\right|\right)\nonumber\\
        &\leq
        \label{equa}
        \epsilon/2 + \sum_{k=1}^m\tilde{\pi}_k\int \left|\prod_{j=1}^d f_{kj}-\prod_{j=1}^d\tilde{f}_{kj}\right|.
    \end{align}
    For all $k\in[m]$, we have the following relation:
    \begin{align*}
        \int\left|\prod_{j=1}^d f_{kj}-\prod_{j=1}^d \tilde{f}_{kj}\right|&\leq \int f_{k1}\left|\prod_{j=2}^df_{kj}-\prod_{j=2}^d \tilde{f}_{kj}\right|+\int \left|f_{k1}-\tilde{f}_{k1}\right|\prod_{j=2}^d \tilde{f}_{kj}\\
        &\leq \epsilon/2d + \int \left|\prod_{j=2}^d f_{kj}-\prod_{j=2}^d\tilde{f}_{kj}\right|\\
        &\leq \epsilon/d + \int\left|\prod_{j=3}^d f_{kj}-\prod_{j=3}^d\tilde{f}_{kj}\right|\leq \dots\leq\epsilon/2.
    \end{align*}
    Combining this with \eqref{equa}, we have $\|f-\tilde{f}\|_1\leq \epsilon$. Thus $\mathcal{N}$ is an $\epsilon$-covering of $\mathcal{G}_{L,q}^{(m,d)}$. Now we calculate the cardinality of $\mathcal{N}$:
    \begin{equation}
\label{number}
|\mathcal{N}| = |S|^{dm}\cdot |D_{\epsilon/2}|\leq|S|^{dm}\cdot \left(\frac{10m}{\epsilon}\right)^{m-1}.
    \end{equation}
    The inequality is from the classical result about the covering number of a simplex (see e.g., Lemma A.4 in \cite{ghosal_entropies_2001}). From the entropic bound of 1-dimensional H\"{o}lder class, we have \cite[see e.g.,][Theorem 27.14]{Polyanskiy_Wu_2025}
     \begin{equation}
         \label{S_upper}
         \log |S| \asymp_{L,q} \left(\frac{d}{\epsilon}\right)^{1/q}.
     \end{equation}
     Thus, plug  \eqref{S_upper} into \eqref{number}  we have
     $$\log N(\mathcal{G}_{L,q}^{(m,d)},\textnormal{TV},\epsilon)\lesssim_{L,q} md\left(\frac{d}{\epsilon}\right)^{1/q}+ (m-1)\log \frac{10m}{\epsilon}\lesssim_{L,q} md^{1+\frac{1}{q}}\left(\frac{1}{\epsilon}\right)^{1/q},$$
     which proves the $\TV$ upper bound.
\par  Now we prove the upper bound under $H$. The idea of choosing a covering set is similar. Pick an $\epsilon/\sqrt{d}$-covering of  $\mathcal{F}_{L,q}$ under $H$, and an $\epsilon^2$-covering of $\Delta^{m-1}$ under $\TV$, denoted by $\mathcal{N}_{\epsilon/\sqrt{d},H},D_{\epsilon^2}$. 
    The covering set is defined as
     $$\mathcal{N}_1=\left\{\tilde{f}=\sum_{k=1}^m\tilde{\pi}_k\prod_{j=1}^d\tilde{f}_{kj}(x_j):\tilde{f}_{kj}\in \mathcal{N}_{\epsilon/\sqrt{d},H}, (\tilde{\pi}_1,\dots,\tilde{\pi}_m)\in D_{\epsilon^2}\right\}.$$
    Now we prove $\mathcal{N}_1$ is an $\epsilon$-covering. For $f\in\mathcal{G}_{\mathcal{F}}$, we pick the element in $\mathcal{N}$ such that
        $$H(f_{kj},\tilde{f}_{kj})\leq \epsilon,\quad\|\pi-\tilde{\pi}\|_1\leq \epsilon^2.$$
        Then we can upper bound $H^2(f,\tilde{f})$:
        \begin{align*}
            H^2(f,\tilde{f})&\leq \left(H\left(\sum_{k=1}^m\pi_k\prod_{j=1}^d f_{kj},\sum_{k=1}^m\pi_k\prod_{j=1}^d\tilde{f}_{kj}\right)+H\left(\sum_{k=1}\pi_k\prod_{j=1}^d\tilde{f}_{kj},\sum_{k=1}^m\tilde{\pi}_k\prod_{j=1}^d\tilde{f}_{kj}\right)\right)^2\\
           & \leq 2\left(H^2\left(\sum_{k=1}^m\pi_k\prod_{j=1}^d f_{kj},\sum_{k=1}\pi_k\prod_{j=1}^d\tilde{f}_{kj}\right)+H^2\left(\sum_{k=1}^m\pi_k\prod_{j=1}^d\tilde{f}_{kj},\sum_{k=1}^m\tilde{\pi}_k\prod_{j=1}^d\tilde{f}_{kj}\right) \right)\\
           &\leq 2\left( \sum_{k=1}^m\pi_kH^2\left(\prod_{j=1}^d f_{kj},\prod_{j=1}^{d}\tilde{f}_{kj}\right)+2\textnormal{TV}\left(\sum_{k=1}^m\pi_k\prod_{j=1}^d\tilde{f}_{kj},\sum_{k=1}^m\tilde{\pi}_k\prod_{j=1}^d\tilde{f}_{kj}\right)\right)\\
           &\leq  2\sum_{k=1}^m \pi_kH^2(\prod_{j=1}^d f_{kj},\prod_{j=1}^{d}\tilde{f}_{kj})+ 4\|\pi-\tilde{\pi}\|_1. 
        \end{align*}
        The first inequality uses the triangle inequality of $H$ as a distance, the second uses the Cauchy-Schwarz inequality, the third uses the convexity of Hellinger distance, and $\frac{H^2}{2}\leq \textnormal{TV}$. Now we bound $H^2(\prod_{j=1}^d f_{kj},\prod_{j=1}^{d}\tilde{f}_{kj})$:
        \begin{align*}
            H^2\left(\prod_{j=1}^d f_{kj},\prod_{j=1}^{d}\tilde{f}_{kj}\right) &= 2\left(1-\prod_{j=1}^d(1-\frac{H^2(f_{kj},\tilde{f}_{kj})}{2})\right)\\
            &\leq 2\left(1-(1-\frac{\epsilon^2}{2d})^d\right)\leq 2\left(1-(1-\frac{\epsilon^2}{2})\right)=\epsilon^2 .
        \end{align*}
        The last inequality is due to $(1+x)^n>nx$ for $x>-1$. Thus, $H^2(f,\tilde{f})\lesssim \epsilon^2$. 
        \par Now we calculate the cardinality of $\mathcal{N}_1$. Similar to \eqref{number}, we have
        \begin{equation}
            \label{N_card}
            |\mathcal{N}_1| = |\mathcal{N}_{\epsilon/\sqrt{d},\mathcal{F}}|^{md}|\mathcal{N}_{\epsilon^2,\Delta^{m-1}}|\lesssim |\mathcal{N}_{\epsilon/\sqrt{d},\mathcal{F}}|^{md}\left(\frac{1}{\epsilon^2}\right)^{m-1}.
        \end{equation}
        Moreover, $\log |\mathcal{N}_{\epsilon/\sqrt{d},\mathcal{F}}|$ has an upper bound given by the entropic bounds of H\"{o}lder class \cite[see][equation (32.56)]{Polyanskiy_Wu_2025}:
        \begin{equation}
            \label{H_N}
            \log |\mathcal{N}_{\epsilon/\sqrt{d},\mathcal{F}}|\asymp_{L,q} d^{1/q}\left(\frac{1}{\epsilon}\right)^{2/q}
        \end{equation}
        Thus, plug \eqref{H_N} into \eqref{N_card} we have the entropic upper bound
        $$\log |\mathcal{N}_1|= md\left(\frac{1}{\epsilon/\sqrt{d}}\right)^{2/q} + (m-1)\log\frac{1}{\epsilon^2}+c\lesssim md^{1+1/q}\left(\frac{1}{\epsilon}\right)^{2/q}$$
        as desired.
        
\noindent\textbf{Lower bound:} We first prove the lower bound under $\TV$. 
For $k=1,2,\dots ,m$, let $\mathcal{F}_{k}$ be a subset of $\mathcal{F}_{L,q}$ such that
\begin{equation}
        \label{sep_supp}
        \mathcal{F}_k\triangleq\left\{f\in\mathcal{F}_{L,q}:\textnormal{supp}(f)\subset \left[\frac{k-1}{m},\frac{k}{m}\right]\right\}.
\end{equation}
For every $k\in[m]$, pick a $2\epsilon$-packing of $\mathcal{F}_{k}$, denoted by $\mathcal{M}_k$.
We consider the following class:
$$P_k\triangleq\left\{p(x)=\prod_{j=1}^d h_{j}(x_j):h_{j}\in \mathcal{M}_k\right\}.$$
We write $P_k=\{p_{1}^{(k)},\dots ,p_{|P_k|}^{(k)}\}$ for $k=1,2,\dots ,m$. Let $M_0=\min_{k\in[m]} |P_k|$, and $A_{M_0,m} =\{1,2,...,M_0\}^m$. 
Now we consider the following packing set:
$$\mathcal{M}:=\left\{g= \frac{1}{m}\sum_{k=1}^m p_{i_k}^{(k)}: p_{i_k}^{(k)}\in P_k,(i_1,\dots ,i_m)\in P_{M_0,m}(\lceil m/2\rceil)\right\},$$
where $P_{M_0,m}(\lceil m/2\rceil)$ is defined as in Proposition \ref{GV_bound}. 
Clearly $\mathcal{M}\subset \mathcal{G}_{\calF_{L,q}}^{(m,d)}$.   

Now we show that $\mathcal{M}$ is an $\epsilon/2$-packing of $\mathcal{G}_{\calF_{L,q}}^{(m,d)}$. 
We first consider the lower bound of $\TV (p_{i}^{(k)},p_{i'}^{(k)})$  for $p_{i}^{(k)},p_{i'}^{(k)}\in P_k,i\neq i'$. Let $p_{i}^{(k)} = \prod_{j=1}^d h_{j}^{(i)},p_{i'}^{(k)} = \prod_{j=1}^d h_{j}^{(i')}$. Since $p_{i}^{(k)}\neq$ there exists a $j_0\in[d]$ such that $h_{j_0}^{(i)}\ne h_{j_0}^{(i')}$. 
Without loss of generality, take $j_0=1$. 
We obtain
\begin{align}
         \TV(p_i^{(k)},p_{i'}^{(k)}) &= \frac{1}{2}\left\|\prod_{j=1}^dh_{j}^{(i)} - \prod_{j=1}^dh_{j}^{(i')} \right\|_1 = \frac{1}{2}\int \left|\prod_{j=1}^dh_{j}^{(i)}(x_j) - \prod_{j=1}^dh_{j}^{(i')}(x_j)\right|\diff x_1\dots \diff x_d\nonumber\\
         &\geq \frac{1}{2}\int \left|\int \left(\prod_{j=1}^dh_{j}^{(i)}(x_j) - \prod_{j=1}^dh_{j}^{(i')}(x_j)\right)\diff x_2\dots \diff x_d\right|\diff x_1\nonumber\\
         \label{TV_p_lower}
         &=\frac{1}{2}\int |h_1^{(i)}(x_1)-h_{1}^{(i')}(x_1)|\diff x_1 =\|h_{1}^{(i)}-h_{1}^{(i')}\|_1\geq \epsilon.
\end{align}
The first inequality is from $|\int f(x)dx|\leq \int |f(x)|dx$. 
The second inequality is due to $h_{1}^{(i)},h_{1}^{(i')}$ are different elements in the packing set $\mathcal{M}_k$.
    
For two different elements $g=\frac{1}{m}\sum_{k=1}^m p_{i_k}^{(k)},\ g' = \frac{1}{m}\sum_{k=1}^m p_{i'_k}^{(k)}\in \mathcal{M}$, the index $i_g=(i_1,\dots ,i_m)$ and $i_{g'}=(i'_1,\dots ,i'_m)$ are two distinct elements in $ P_{M_0,m}(\lceil m/2\rceil)$. Thus, there exists $ S \subseteq [m], |S|\geq \lceil m/2\rceil$, such that for $k\in S$, $i_{k}\ne i_k'$. 
This implies $p_{i_k}^{(k)}\ne p_{i_k'}^{(k)}$. 
From \eqref{TV_p_lower}, we deduce that
    \begin{align*}
       \textnormal{TV}(g,g') &= \frac{1}{2m}\left\|\sum_{k=1}^m(p_{i_k}^{(k)}-p_{i'_k}^{(k)})\right\|_1 \\
       &=\frac{1}{2m}\sum_{k=1}^m\left\|p_{i_k}^{(k)}-p_{i'_k}^{(k)}\right\|_1 \quad \text{(the support of components are disjoint)}\\
       & \geq \frac{1}{2m}\sum_{k\in S}\left\|p_{i_k}^{(k)}-p_{i'_k}^{(k)}\right\|_1 \geq \frac{1}{2m}\cdot 2\epsilon\lceil m/2\rceil\geq \epsilon/2.
    \end{align*}
Hence, $\mathcal{M}$ is an $\epsilon/2$-packing of $\mathcal{G}_{\calF_{L,q}}^{(m,d)}$. 
Now we calculate the cardinality of $\mathcal{M}$, given by $|\mathcal{M}| = | P_{M_0,m}(\lceil m/2\rceil)|$.
Applying Proposition \ref{GV_bound}, we obtain
    $$\left| P_{M_0,m}(\lceil m/2\rceil)\right| \geq \frac{M_0^m}{\sum_{j=0}^{\lceil m/2\rceil-1}\binom{m}{j}(M_0-1)^j}.$$
    Applying the inequality $\binom{m}{j}\leq \binom{m}{\lfloor \frac{m}{2}\rfloor}\leq (\frac{me}{m/2})^{m/2}\leq (\sqrt{2e})^m$, we have
\begin{equation}
        \label{combination_bound}
        | P_{M_0,m}(\lceil m/2\rceil)|\geq \frac{M_0^m}{(\sqrt{2e})^m \sum_{j=0}^{\lceil m/2\rceil-1}(M-1)^j}\geq \frac{M_0^m}{(\sqrt{2e})^m M_0^\frac{m}{2}}\geq \left(\sqrt{\frac{M_0}{2e}}\right)^m.
\end{equation}
We have the following lower bound for $M_0$:
    \begin{equation}
        \label{M_0_lower}
        \log M_0 = \min_{k\in[m]} \log|P_k|= \min_{k\in[m]} d\log|\mathcal{M}_k| \gtrsim_{L,q} \frac{d}{m}\left(\frac{1}{\epsilon}\right)^{1/q}.
    \end{equation}
The last inequality is from the entropic bound of 1-dimensional H\"{o}lder class \cite[see e.g.,][Theorem 27.14]{Polyanskiy_Wu_2025}.
    Plugging \eqref{M_0_lower} into \eqref{combination_bound}, we have
    $$\log|\mathcal{M}| \gtrsim m\log|M_0|\gtrsim_{L,q} d\cdot\left(\frac{1}{\epsilon}\right)^{1/q}.$$
This completes the proof of the entropic lower bound for $\TV$.

Now we turn to the lower bound for $H$. We pick an $\epsilon/\sqrt{d}$-packing of $\mathcal{F}_{k}$ in \eqref{sep_supp}, denoted by $\mathcal{M}_{k,\epsilon/\sqrt{d}}=\{g_{1}^{(k)},\dots g_{|\mathcal{M}_{k,\epsilon/\sqrt{d}}|}^{(k)}\}$. 
Let $M_1:=\min_{k\in[m]}|\mathcal{M}_{k,\epsilon/\sqrt{d}}|$. We consider the following set:
    $$Q_k = \left\{\prod_{j=1}^d g_{i_j}^{(k)}(x_j):(i_1,\dots ,i_d)\in P_{M_1,d}(\lceil d/2\rceil)\right\},$$
where $P_{M_1,d}(\lceil d/2\rceil)$ is defined as in Proposition \ref{GV_bound}. 
We write $Q_k=\{q_{1}^{(k)},...,q_{|Q_k|^{(k)}}\}$ and let $M_2=\min_{k\in[m]}|Q_k|$. 
Now we construct the packing set to be
    \begin{equation}
        \mathcal{M}_1 = \left\{g=\frac{1}{m}\sum_{k=1}^mq_{i_k}^{(k)}: q_{i_k}^{(k)}\in Q_k,\ (i_1,\dots i_m)\in P_{M_2,m}(\lceil m/2\rceil)\right\}.
    \end{equation}
We now prove $\mathcal{M}_1$ is an $\epsilon/\sqrt{8}$-packing.  For two different elements $q_{i}^{(k)}=\prod_{j=1}^d g^{(k)}_{i_j},q_{i'}^{(k)}=\prod_{j=1}^dg_{i'_j}^{(k)}\in Q_k$, there exists $T\subseteq [d], |T|\geq \lceil d/2\rceil$, such that for all $j\in T$, $i_j\ne i_j'$. 
Thus, we have the lower bound of Hellinger distance:
    \begin{align}
        H^2(q_{i}^{(k)},q_{i'}^{(k)}) &= H^2\left(\prod_{j=1}^d g^{(k)}_{i_j},\prod_{j=1}^dg_{i'_j}^{(k)}\right)=2\left(1-\prod_{j=1}^d\left(1-\frac{H^2(g_{i_j}^{(k)},g^{(k)}_{i_j'})}{2}\right) \right)\nonumber\\
        &\geq 2\left(1-(1-\frac{\epsilon^2}{2d})^{\lceil d/2\rceil}\right)\nonumber\\
        &\geq 2\left( 1- \left(1-\frac{\epsilon^2}{2d}\cdot \lceil d/2\rceil +\binom{\lceil d/2\rceil}{2}\left(\frac{\epsilon^2}{2d}\right)^2\right) \right)\nonumber\\ \label{re_product_dist_H}
        &\geq \frac{\epsilon^2}{2} - \frac{(d+2)d}{4}\cdot \frac{\epsilon^4}{4d^2}
        \geq \epsilon^2/4.
    \end{align}
For two distinct elements $g=\frac{1}{m}\sum_{k=1}^m q_{i_k}^{(k)},g'=\frac{1}{m}\sum_{k=1}^m q_{i_k'}^{(k)}\in\mathcal{M}_1$, there exists $T'\subseteq [m], |T'|\geq \lceil m/2\rceil$, such that for all $k\in T'$, $i_k\ne i_k'$. From \eqref{re_product_dist_H} we deduce that
\begin{align}
    H^2(g,g') &= \frac{1}{m} H^2\left(\sum_{k=1}^mq_{i_k}^{(k)},\sum_{k=1}^mq_{i_k'}^{(k)}\right)\nonumber\\
    &= \frac{1}{m}\sum_{k=1}^m H^2\left(q_{i_k}^{(k)},q_{i_k'}^{(k)}\right)\quad \text{(the support of components are disjoint)}\nonumber\\
    \label{H_inequ}
    &\geq \frac{1}{m}\sum_{k\in S}H^2\left(q_{i_k}^{(k)},q_{i_k'}^{(k)}\right)\geq \frac{1}{m}\cdot\lceil m/2\rceil \frac{\epsilon^2}{4}\geq \epsilon^2/8.
\end{align}
    This proves that $\mathcal{M}_1$ is $\epsilon/\sqrt{8}$-packing. Now we calculate the cardinality of $\calM_1$, given by
    $|\mathcal{M}_1| = \left|P_{M_2,m}(\lceil m/2\rceil)\right|.$
     Similar to \eqref{combination_bound},  we have
    $$|P_{M_2,m}(\lceil m/2\rceil)|\geq \left(\sqrt{\frac{M_2}{2e}}\right)^m,\quad M_2 =|P_{M_1,d}(\lceil d/2\rceil)|\geq \left(\sqrt{\frac{M_1}{2e}}\right)^d.$$
    This implies
    \begin{equation}
        \label{M_1_lower_bound}
        \log|\mathcal{M}_1|\gtrsim md\log M_1 \gtrsim_{L,q} md\cdot\frac{1}{m}\cdot d^{1/q}\left(\frac{1}{\epsilon}\right)^{2/q}=d^{1+\frac{1}{q}}\left(\frac{1}{\epsilon}\right)^{2/q}.
    \end{equation}
    The last inequality is given by the entropic bounds of H\"{o}lder class \cite[see][equation (32.56)]{Polyanskiy_Wu_2025}.
\end{proof}

\subsection{Wrapping up the proof}
\noindent With the entropic bounds in Lemma \ref{entro_lower_bound}, we are ready to prove Theorem \ref{Holder_rate}.
\begin{proof}[Proof of Theorem \ref{Holder_rate}]
    The upper bound in Theorem \ref{Holder_rate} is directly from Proposition \ref{results_mini_upper_classic}. Let $V_{\rho}(\epsilon)$ be an upper bound of the $\epsilon$-covering entropy under $\rho$, for $\rho=H$, we have 
$$R^*_{H} \lesssim \inf_{\epsilon>0}\left\{ \epsilon^2+\frac{1}{n}V_{H}(\epsilon)\right\}\lesssim_{L,q} \inf_{0<\epsilon<1}\left\{ \epsilon^2 +  md^{1+\frac{1}{q}}\left(\frac{1}{\epsilon}\right)^{2/q}\right\}$$
    Let $\epsilon=\epsilon_{n,m,d} = n^{-\frac{q}{2q+2}}m^{\frac{q}{2q+2}}d^{\frac{1}{2}}$ to get the minimax upper bound for $H$. To guarantee $\epsilon_{n,m,d}<1$, we need $n\geq md^{1+\frac{1}{q}}$. Similarly, we can derive the minimax upper bound for $\TV$. We omit the details here.
\begin{equation}
    \label{subclass}
         \mathcal{G}_{L,q,1}^{(m,d)}:= \left\{f = \prod_{j=1}^d f_j: f_j\in\mathcal{F}_{L,q}\right\}.
\end{equation}
Then, $\mathcal{G}_{L,q,1}^{(m,d)}\subset \mathcal{G}_{\calF_{L,q}}^{(m,d)}$ and thus
$$R^*_{H,\calF_{L,q}}\geq \inf_{\hat{f}_n}\sup_{f\in\mathcal{G}_{L,q,1}^{(m,d)}}\Expect\left[H^2(\hat{f}_n,f)\right].$$
We will calculate the covering radius of $\calG_{L,q,1}^{(m,d)}$ defined in Proposition \ref{Yang_minimax}.
Now we pick an $\KL$ $\epsilon/\sqrt{d}$-covering of $\mathcal{F}_{L,q}$, denoted by $\mathcal{N}_{\KL}$.
     We consider the following set:
     $$\mathcal{N}_2 = \left\{\tilde{f} = \prod_{j=1}^d \tilde{f}_j(x_j): \tilde{f}_j\in \mathcal{N}_{\KL}\right\}.$$
     We will now show that $\mathcal{N}_{2}$ is an $\KL$ $\epsilon$-covering of $\mathcal{G}_{L,q,1}^{(m,d)}$. For any $f\in\mathcal{G}_{L,q,1}^{(m,d)}$, we find an element $\tilde{f}\in \mathcal{N}_2$ such that $\sqrt{\KL(f_j,\tilde{f}_j)}\leq \epsilon/\sqrt{d}$.
     From the additivity of $\KL$-divergence for product density, we have
$$\sqrt{\KL(f,\tilde{f})} = \sqrt{\sum_{j=1}^d \KL(f_j,\tilde{f}_j)} \leq \epsilon.$$
   This shows that
    \begin{equation}
    \label{N_1}
        \log N(\mathcal{G}_{L,q,1}^{(m,d)},\sqrt{\KL},\epsilon)\leq d\log N(\mathcal{F}_{L,q}, \sqrt{\KL}, \epsilon/\sqrt{d}).
    \end{equation}
     Now we derive an upper bound for $\KL$ covering entropy of $\mathcal{F}_{L,q}$. We claim that the density class $\mathcal{F}_{L,q}$ has a finite $\chi^2$ radius:
     \begin{equation*}
         \inf_{u}\sup_{f\in\mathcal{F}_{L,q}} \chi^2(f||u)<\infty.
     \end{equation*}
     This can be verified by choosing $u$ the density of uniform distribution on $[0,1]$:
    \begin{equation*}
        \inf_{u}\sup_{f\in\mathcal{F}_{L,q}} \chi^2(f||u)\leq \sup_{f\in\mathcal{F}_{L,q},u\sim \textnormal{Unif}[0,1]} \chi^2(f||u)=\sup_{f\in\mathcal{F}_{L,q}} \int f(x)^2dx -1<\infty.
    \end{equation*}
     Thus, by  Theorem 32.6 in  \cite{Polyanskiy_Wu_2025} with $\lambda=2$, we have 
     \begin{equation*}
         N\left(\mathcal{F}_{L,q},\sqrt{\textnormal{\KL}},\epsilon\sqrt{\log\frac{1}{\epsilon}}\right) \lesssim_{L,q} N(\mathcal{F}_{L,q},H,\epsilon).
     \end{equation*}
     Combining this with \eqref{N_1}, we have  
     $$ \log N(\mathcal{G}_{L,q,1}^{(m,d)},\sqrt{\KL},\epsilon)\leq d\log N(\mathcal{F}_{L,q}, \sqrt{\KL}, \epsilon/\sqrt{d})\lesssim_{L,q} d \log N(\mathcal{F}_{L,q},H, \delta/\sqrt{d}):= V_H(\delta),$$
     where $\delta$ satisfies $\epsilon=\delta\sqrt{\log \frac{1}{\delta}}$. Now we calculate covering radius of $\calG_{L,q,1}^{(m,d)}$. We know $V_H(\delta_n) \gtrsim_{L,q} n\epsilon_n^2$ for $\epsilon_n=\delta_n\sqrt{\log\frac{1}{\delta_n}}$, thus
     $$d^{1+1/q}\left(\frac{1}{\delta_n}\right)^{2/q}\gtrsim_{L,q} n\delta_n^2\log \frac{1}{\delta_n},$$
which gives $n\epsilon_n^2  \lesssim_{L,q} d(n\log n)^{\frac{1}{q+1}}$.
Now we apply Proposition \ref{Yang_minimax} to obtain the minimax lower bound.
From Lemma \ref{entro_lower_bound} we know
     \begin{equation*}
         \log M(\mathcal{G}_{L,q,1}^{(m,d)},H,\epsilon)\gtrsim_{L,q} d^{1+\frac{1}{q}}\left(\frac{1}{\epsilon}\right)^{2/q}.
     \end{equation*}
      Now, plug this and the formula of $n\epsilon_n^2$ into \eqref{Packing_lower}, we have 
     \begin{equation*}
         \label{epsilon_nH}
         d^{1+\frac{1}{q}}\left(\frac{1}{\epsilon_{n,H}}\right)^{2/q}\lesssim_{L,q}d(n\log n)^{\frac{1}{q+1}}\implies\epsilon_{n,H}^2 \gtrsim_{L,q} d(n\log n)^{-\frac{q}{q+1}}.
     \end{equation*}
This proves the minimax lower bound under $H$. 
For $\TV$, from Lemma \ref{entro_lower_bound} again, 
$$\log M(\mathcal{G}_{L,q,1}^{(m,d)},\TV,\epsilon)\gtrsim_{L,q} d\left(\frac{1}{\epsilon}\right)^{1/q}.$$
Thus,
$$d\left(\frac{1}{\epsilon_{n,\TV}}\right)^{1/q}\lesssim_{L,q}d(n\log n)^{\frac{1}{q+1}}\Rightarrow\epsilon_{n,\TV}^2\gtrsim_{L,q}(n\log n)^{-\frac{q}{q+1}}.$$
\end{proof}

\section{Details in Section \ref{sec_algo}}
\subsection{Recovering the component from the exact joint density}
\label{algo_compare}
 \par In this subsection, we present the recovery procedure from the known joint density $f$ and discuss its connection to Algorithm 1. The joint density can be expressed as
 \begin{equation}
     \label{joint_blocked}
     f(x_1,\dots ,x_{2m-1}) = \sum_{k=1}^m\pi_k f_{k}^{(1)}(y)f_{k}^{(2)}(z)f_{k(2m-1)}(x_{2m-1}),
 \end{equation}
 where $y=(x_1,\dots ,x_{m-1}),z = (x_{m},\dots ,x_{2m-2})$ and $f_{k}^{(1)}(y) = \prod_{j=1}^{m-1}f_{kj}(x_j),f_{k}^{(2)}(z) = \prod_{j=m}^{2m-2}f_{kj}(x_j)$.
 Integrating  over $x_{2m-1}$, we obtain
 \begin{equation}
     \label{T_+yz}
     T_{+}(y,z) \triangleq \sum_{k=1}^m\pi_kf_{k}^{(1)}(y) f_{k}^{(2)}(z).
 \end{equation}
  Applying a unitary transformation $\tilde{U}$, we map $T_+(y,z)\in L^2(\real^{m-1}\times \real^{m-1})$ to the following linear operator: 
 \begin{equation}
     \label{T_+_repre}
     T_{+} \triangleq \tilde{U}^{-1}(T_{+}(y,z)) = F_1 D_{\pi} F_2^{*}\in \calB\left(L^2(\real^{m-1}),L^2(\real^{m-1})\right),
 \end{equation}
 where $F_1=(f_{1}^{(1)},\dots ,f_{m}^{(1)}), F_2=(f_{1}^{(2)},\dots ,f_{m}^{(2)})$ and $D_{\pi}=\textnormal{diag}(\pi_1,\dots ,\pi_m)$.
 Since $T_{+}$ is a finite rank operator, we can perform its singular value decomposition (SVD):
 \begin{equation}
     \label{SVD_T}
     T_{+} = \sum_{k=1}^m \lambda_k \phi_k\otimes \psi_{k} = U \Sigma V^*,
 \end{equation}
 where $U=(\phi_1,\dots ,\phi_{m}),V=(\psi_1,\dots ,\psi_m)$ are orthonormal and $\Sigma = \diag(\lambda_1,\dots ,\lambda_m)$.
 Since $\{f_{kj}\}_{j=1}^{2m-1}$ are $\mu$-incoherent, hence pairwise distinct, $F_1$ and $F_2$ both have full column rank, implying that the diagonal entries of $\Sigma$ are positive.
 Let $T_{+}^{\dagger}$ denote the Moore-Penrose inverse of $T_{+}$, given explicitly by
  \begin{equation}
     \label{T_inverse}
     T_{+}^{\dagger} = (F_{2}^*)^{\dagger} D_{\pi}^{-1} F_1^{\dagger} =V \Sigma^{-1}U^*.
 \end{equation}
 We now select a subset $A$ of the support of the $(2m-1)$-th variable and define the operator
  \begin{equation}
      \label{T_A}
      T_A \triangleq \tilde{U}^{-1}\left(\int_A f(x_1,\dots,x_{2m-1})dx_{2m-1}\right) = F_1D_{\pi,A} F_2^*,
  \end{equation}
   where $D_{\pi,A} = \diag(\pi_1a_1,\cdots,\pi_ma_m)$ with $a_k=\int_A f_{k(2m-1)}(x)dx$ for $k=1,2,\dots,m$.
We have the following result.
\begin{lemma}
\label{lemma:simu_diag}
    Let $T_+,T_A$ be defined as in \eqref{T_+_repre},\eqref{T_A}, respectively. Then for each $k=1,2,\cdots,m$, $f_{k}^{(1)}$ is eigenfunction of $T_AT^{\dagger}_+$. Moreover, if $a_k=\int_A f_{k(2m-1)}(x)dx$ are pairwise distinct for $k=1,2,\dots,m$, then up to a permutation, $T_AT_+$ uniquely determines $f_{k}^{(1)}$. 
\end{lemma}
\begin{proof}
      From \eqref{T_inverse},\eqref{T_A}, we have $ T_AT_{+}^{\dagger}F_1 =  F_1 D_{\pi,A} (F_{2}^*)^{\dagger} D_{\pi}^{-1} F_1^{\dagger}F_1 = F_1 \diag(a_1,\dots,a_m)$. If $a_k$'s are pairwise distinct, then the eigenspaces are one-dimensional, and each $f_{k}^{(1)}$ is determined uniquely up to scaling. Since $f_{k}^{(1)}$ is a density function, the normalization further fixes it.
\end{proof}
Lemma \ref{lemma:simu_diag} shows that  $f_{k}^{(1)}$'s are eigenfunctions of $T_AT_+^{\dagger}$.
Consequently, 
$F_1$ simultaneously diagonalizes $T_AT_+^{\dagger}$ for any choice of $A$.
In practice, instead of working directly with 
$T_AT_+$, we compute its coefficient matrix under the basis 
$U$:
\begin{equation}
    \label{eta_A}
    \eta_A\triangleq U^*T_AT^{\dagger}_{+}U \in \real^{m\times m}.
\end{equation}
Let $W$ be the matrix whose columns are the eigenvectors of $\eta_A$. Then $W$ represents the coefficients of $F_1$ under the basis $U$. Thus,
\begin{equation}
    \label{F_1_formula}
    (g_{1},\dots ,g_{m}) = UW, \quad F_1 = (f_{1}^{(1)},\dots ,f_{m}^{(1)}) = (g_{1}/\|g_{1}\|_1,\dots ,g_{m}/\|g_{m }\|_1).
\end{equation}
We summarize the above procedure in Algorithm 2 below.
 \begin{algorithm}
    \caption{Recover the component density from true density}
    \begin{algorithmic}[1]
    \Require Joint density $f$ that admits model \eqref{MM_den}.
    \Ensure $F_1=(f_{1}^{(1)},\dots ,f_{m}^{(1)})$ 
    \State  Calculate $T_{+}(y,z) = \int f(y,z,x_{2m-1})\diff x_{2m-1}$, where $y=(x_1,\dots ,x_{m-1})$ and $z=(x_m,\dots ,x_{2m-2})$.   
    \State Perform SVD on $T_{+} = U\Sigma V^*$. Let $T_+^{\dagger} = V\Sigma^{-1}U^*$
     \State Choose some subset $A$, let $T_A = \int_A f(y,z,x_{2m-1})\diff x_{2m-1}$
    \State  Let $\eta_A=U^*T_AT_{+}^{\dagger}U$, calculate $W=(w_1,\dots ,w_m)$ the columns of $L^2$ unit eigenvectors of $\eta_A$
     \State Let $(g_{1},\dots ,g_{m}) = UW$, return $ F_1 = (f_{1}^{(1)},\dots ,f_{m}^{(1)}) = (g_{1}/\|g_{1}\|_1,\dots ,g_{m}/\|g_{m }\|_1)$
    \end{algorithmic}
\end{algorithm}  
Finally, note that Algorithm 1 in the main text is simply a plug-in version of Algorithm 2.
\subsection{Proof of Theorem \ref{algo_error}}
\label{proof_algo}
We need the following perturbation lemmas. The first one is for eigenvectors, and the second is for pseudo pseudo-inverse of linear operators.
\begin{lemma}[Theorem 2.8 in \cite{stewart1990matrix}]
\label{eigenvec_error}
    Let $A$ be a diagonalizable real matrix with eigen decomposition $U^{-1}AU=D$. Rewrite the decomposition as follows:
    $$(v_1,V_2)^{*} A (u_1,U_2) = \begin{pmatrix}
        \lambda_1 &0\\0 &L_2
    \end{pmatrix},$$
    where $U=(u_1,U_2)$ and $(v_1,V_2)^* = (u_1,U_2)^{-1}$. Then for $\hat{A}=A+E, \|E\|\leq \epsilon$, we have
    \begin{equation}
        \label{eigenvec_bound}
        \|u_{1}-\hat{u}_1\|\leq C_1 \|U_2(\lambda_1 I-L_2)^{-1}V_2^{*}\|\epsilon,
    \end{equation}
    where $C_1>0$ is an absolute constant.
\end{lemma}

\begin{lemma}[Theorem 2 in \cite{chen_expression_1998}]
\label{ginv_bound}
    Let $\mathcal{H}_1,\mathcal{H}_2$ be Hilbert spaces and $T,\hat{T}$ be two linear operators from $\mathcal{H}_1$ to $\mathcal{H}_2$. Suppose $\hat{T}=T+E$ such that $\textnormal{rank}(T) =\textnormal{rank}(\hat{T})<\infty$. Then 
    \begin{equation}
        \label{p_inverse_bound}
        \frac{\|\hat{T}^{\dagger}-T^{\dagger}\|}{\|T^{\dagger}\|}\leq \frac{3\|T^{\dagger}\|\|E\|}{1-\|\hat{T}^{\dagger}\|\|E\|}.
    \end{equation}
\end{lemma}

\begin{proof}[Proof of Theorem \ref{algo_error}]
    In this proof, the norm $\|\cdot\|$ refers to the operator norm if not specified. The notation, if not followed by the name of the variable, should be understood as elements in the tensor of Hilbert spaces, like the relationship between $T_{+}(y,z)$ and $T_{+}$ in \eqref{T_+yz}, \eqref{T_+_repre}. The operator norm in the tensor of Hilbert spaces is identical to the $L^2$ norm in the $L^2$ function space, because the transform between the two spaces is unitary. 
    \par We write $\hat{f}(y,z,x_{2m-1})  = f(y,z,x_{2m-1})+E(y,z,x_{2m-1})$, such that
    \begin{equation}
        \label{decom}
        f(y,z,x_{2m-1})=\sum_{k=1}^{m}\pi_k f_{k}^{(1)}(y)f_{k}^{(2)}(z)f_{k(2m-1)}(x_{2m-1}), \|E(y,z,x_{2m-1})\|_2\leq \epsilon.
    \end{equation}
    Let $\hat{T}_{+,m} = \sum_{k=1}^{m}\hat{\lambda}_k\hat{\phi}_{k}\otimes \hat{\psi}_k$. We can obtain that $\hat{T}_{+,m}$ is close to $T_{+}$ in \eqref{T_+_repre}:
    \begin{equation}
        \label{op_rela}
        \|\hat{T}_{+,m}-T_{+}\| \leq \|\hat{T}_{+,m}-\hat{T}_{+}\|  + \|\hat{T}_+ - T_+\|\leq 2\|\hat{T}_{+}-T_+\| . 
    \end{equation}
    The first inequality is due to the triangle inequality, the second due to the choice of $\hat{T}_{+,m}$ and the fact that $T_{+}$ is rank $m$. Now, from Cauchy-Schwarz inequality, we bound the right-hand side of \eqref{op_rela}:
    \begin{equation}
        \label{upper_E3}
        \|\hat{T}_{+} - T_{+}\| = \|\hat{T}_{+}(y,z)-T_{+}(y,z)\|_2=\sqrt{\int\int\left(\int E(y,z,x_{2m-1})dx_{2m-1}\right)^2 dydz}\leq \|E\|_2\leq\epsilon.
    \end{equation}
    Thus $\|\hat{T}_{+,m} - T_+\| \leq 2\epsilon$. The $m$-th singular value of $T_{+}$ is lower bounded from the equation \eqref{T_+_repre} and \eqref{condition_num}:
    \begin{equation}
        \label{sigma_lower}
        \sigma:=\sigma_{m}(T_{+})\geq \sigma_m(F_1)\sigma_{m}(D_{\pi})\sigma_m(F_2^*)\geq \frac{\zeta (1-\mu)^m}{(m-1)!}.
    \end{equation}
    From the condition in Theorem \ref{algo_error}, we have $\sigma\geq 4\epsilon$. Thus, from Lemma \ref{Weyl} we have
    \begin{equation}
        \label{hat_T_singular}
        |\sigma-\sigma_m(\hat{T}_{+,m})|\leq 2\epsilon\Rightarrow \sigma_m(\hat{T}_{+,m})\geq\frac{1}{2}\sigma.
    \end{equation}
    Now we apply Lemma \ref{ginv_bound} to obtain
    \begin{equation}
        \label{re_bound_inverse}
        \|\hat{T}_{+,m}^{\dagger}-T_{+}^{\dagger}\| \leq \frac{3\|T_{+}^{\dagger}\|^2 \|\hat{T}_{+,m} - T_+\|}{1-\|\hat{T}_{+,m}^{\dagger}\|\|\hat{T}_{+,m} - T_+\|}\leq \frac{3\epsilon/\sigma^2}{1-\frac{2}{\sigma}\cdot \frac{\sigma}{4}} \leq \frac{6\epsilon}{\sigma^2}.
    \end{equation}
    Let $\hat{T}_A(y,z) = \int_A \hat{f}(y,z,x_{2m-1})dx_{2m-1}$. Now we calculate the error between $\hat{T}_A$ and $T_A$ in \eqref{T_A}. From the Cauchy-Schwarz inequality again, we have
    \begin{align}
    \label{upper_A}
        \|T_A-\hat{T}_A\| &= \|T_A(y,z)-\hat{T}_A(y,z)\|_2\nonumber\\
        &= \sqrt{\int(\int_A f(y,z,x_{2m-1})dx_{2m-1}- \int_A \hat{f}(y,z,x_{2m-1})dx_{2m-1})^2dydz}\nonumber\\
        &\leq \sqrt{\frac{1}{\mu_{Leb}(A)}}\|E\|_2\leq \frac{\epsilon}{\sqrt{\mu_0}}.
    \end{align}
    Moreover, $T_A$ is upper bounded by a constant $L_{C,m}^{(0)}$ since all $f_{kj}$ are upper bounded by $C$. From  \eqref{re_bound_inverse} and \eqref{upper_A}, we can now give an error upper bound for the object of eigen decomposition $T_AT_{+}^{\dagger}$:
    \begin{align}
        \|T_AT_{+}^{\dagger}-\hat{T}_A\hat{T}_{+,m}^{\dagger}\| &= \|T_AT_{+}^{\dagger}- T_A\hat{T}_{+,m}^{\dagger}+ T_A\hat{T}_{+,m}^{\dagger} -\hat{T}_A\hat{T}_{+,m}^{\dagger}\|\nonumber\\
        &\leq \|T_A\|\|T_+^{\dagger}-\hat{T}_{+,m}^{\dagger}\| +\|T_A-\hat{T}_A\|\|\hat{T}_{+,m}^{\dagger}\|\nonumber\\
        &\leq \frac{6\epsilon}{\sigma^2}\|T_A\| + \frac{2\epsilon}{\sqrt{\mu_0}\sigma}\leq \frac{(6L^{(0)}_{C,m}+2)\epsilon}{\sigma^2\sqrt{\mu_0}}.
    \end{align}
    Let $\hat{U}=(\hat{\phi}_1,\dots ,\hat{\phi}_m)$, next we need to upper bound the error of $U$ in \eqref{SVD_T} and $\hat{U}$. From \eqref{op_rela}, we have $\|\hat{T}_{+,m} - T_{+}\|\leq 2\epsilon$. Now, from Davis-Kahan Sin$\Theta$ theorem (see e.g., Theorem VII.3.2 in \cite{bhatia_matrix_1997}), we have
    $$\|\sin \left(U,\hat{U}\right)\| \leq \frac{2\epsilon}{\sigma}:=\tilde{\epsilon}_1 \implies \|\cos \left(U,\hat{U}\right)\|= \|U^*\hat{U}\|\geq \sqrt{1-\tilde{\epsilon}_1^2}.$$
    Thus, we have
    \begin{equation}
        \label{error_U}
        \|U-\hat{U}\|=\|U^*U-U^*\hat{U}\|=\|I-U^*\hat{U}\|\leq 1-\sqrt{1-\tilde{\epsilon}_1^2}\leq \tilde{\epsilon}_1.
    \end{equation}
    Now we can upper bound the error between $\eta_A$ in \eqref{eta_A} and $\hat{\eta}_A$:
    \begin{align}
         \label{eta_A_errorbound}
        \|\eta_A-\hat{\eta}_A\| &= \|U^*T_AT_{+}^{\dagger}U - \hat{U}^* T_AT_{+}^{\dagger}U+ \hat{U}^* T_AT_{+}^{\dagger}U - \hat{U}^{*}\hat{T}_A\hat{T}_{+,m} U+\hat{U}^{*}\hat{T}_A\hat{T}_{+,m} U - \hat{U}^{*}\hat{T}_A\hat{T}_{+,m} \hat{U}\|\nonumber\\
        &\leq \|U-\hat{U}\|\|T_AT_{+}^{\dagger}\| + \|T_AT_{+}^{\dagger}-\hat{T}_A\hat{T}_{+,m}^{\dagger}\| + \|U-\hat{U}\|\|\hat{T}_A\hat{T}_{+,m}\|\nonumber\\
        &\leq \tilde{\epsilon}_1(\|T_AT_+^{\dagger}\| +\|\hat{T}_A\hat{T}_{+,m}^{\dagger}\|) + \frac{(6L_{C,m}^{(0)}+2)\epsilon}{\sigma^2\sqrt{\mu_0}}\nonumber\\
        &\leq \frac{L_{C,m}^{(1)}\tilde{\epsilon}_1}{\sqrt{\mu_0}\sigma} + \frac{(6L_{C,m}^{(0)}+2)\epsilon}{\sigma^2\sqrt{\mu_0}}\leq \frac{L_{C,m}^{(2)}\epsilon}{\sigma^2\sqrt{\mu_0}}:=\tilde{\epsilon}_2,
    \end{align}
    Now we are ready to upper bound the error between $\hat{W}$ and $W$ in \eqref{F_1_formula}. In \eqref{F_1_formula}, we know $U$ and $F_1$ are both full column rank, thus $W$ is invertible. We now write the eigen decomposition of true $\eta_A$:
    \begin{equation}
        \label{eigen_eta}
        (v_k,V_{-k})^{*}\eta_A (w_k,W_{-k}) = \begin{pmatrix}
            \lambda_k&0\\ 0& L_{-k}
        \end{pmatrix},
    \end{equation}
    where $W=(w_k,W_{-k})$ and $V^*=(v_k,V_{-k})^*=W^{-1}$. We know $\|W_{-k}\|=1, \|V_{-k}^{H}\|\leq \|V^{H}\| = \|W^{-1}\| \leq 1/\sigma$. Thus, applying Lemma \ref{eigenvec_error} combined with \eqref{eta_A_errorbound} we have 
    \begin{equation}
        \label{w_bound}
        \|w_k-\hat{w}_k\|\leq C_1 \|W_{-k}(\lambda_k I-L_{-k})V_{-k}^{H}\|\tilde{\epsilon}_2 \leq  \frac{L_{C,m}^{(3)}\epsilon}{\sigma^3\delta\sqrt{\mu_0}}
    \end{equation}
    for some constant $C_{m,1}>0$. Now, let $g_{k}$ be the functions in \eqref{F_1_formula}, for a constant $C_2>0$ we have
    \begin{align*}
        \label{upper_f1}
        \|g_{k}-\hat{g}_{k}\|_2 &=\|Uw_k-\hat{U}\hat{w}_k\|\leq \|Uw_k-\hat{U}w_k\|+\|\hat{U}w_k-\hat{U}\hat{w}_k\|\\
        &\leq \|U-\hat{U}\|+\|w_k-\hat{w}_k\| \leq \frac{2\epsilon}{\sigma}+ L_{C,m}^{(3)} \frac{\epsilon}{\sigma^3\delta\sqrt{\mu_0}}\leq \frac{L_{C,m}^{(4)}\epsilon}{\sigma^3\delta\sqrt{\mu_0}}:=\epsilon_3.
    \end{align*}
     The condition of $\epsilon$ in Theorem \ref{algo_error} ensures the condition of Lemma \ref{L_12_error}. Suppose $g_{k}$ is upper bounded by a constant $L_{C,m}^{(5)}$. We apply Lemma \ref{L_12_error} to obtain 
    $$\|\hat{h}_{k}-f_{k}^{(1)}\|_2 \leq 8(L_{C,m}^{(5)})^2\epsilon_3,$$
    where $f_{k}^{(1)}$ is defined in equation \eqref{joint_blocked}. Now, since $f_{k1}$ is on $[0,1]$, we do the integral and apply Cauchy-Schwarz to obtain
    $$\|\hat{f}_{k1}-f_{k1}\|_2 = \left\|\int \hat{h}_{k}dx_{2}\dots dx_{m-1}-\int h_{k}dx_{2}\dots dx_{m-1}\right\|_2\leq \|\hat{h}_{k}-h_{k}\|_2\leq \frac{L_{C,m}\epsilon}{\sigma^3\delta\sqrt{\mu_0}}.$$ 
    Now plug in the lower bound of $\sigma$ in \eqref{sigma_lower} to obtain the result as desired.
\end{proof}

\end{document}